\label{key}
\documentclass[english, 11pt]{article}

\usepackage[T1]{fontenc}
\usepackage[latin9]{inputenc}
\usepackage{amsthm}
\usepackage{amsmath}
\usepackage{amssymb}
\usepackage{color}
\usepackage{esint}
\usepackage{verbatim}
\usepackage{appendix}

\makeatletter
\numberwithin{equation}{section}
\numberwithin{figure}{section}
\theoremstyle{plain}
\newtheorem{thm}{\protect\theoremname}
\theoremstyle{definition}
\newtheorem{defn}[thm]{\protect\definitionname}
\theoremstyle{remark}
\newtheorem{rem}[thm]{\protect\remarkname}
\theoremstyle{plain}
\newtheorem{cor}[thm]{\protect\corollaryname}
\newtheorem{prop}{Proposition}

\@ifundefined{date}{}{\date{}}
\makeatother

\usepackage{babel}
\providecommand{\corollaryname}{Corollary}
\providecommand{\definitionname}{Definition}
\providecommand{\remarkname}{Remark}
\providecommand{\theoremname}{Theorem}

\title{Inviscid incompressible limit \\ for compressible micro-polar fluids}
\author{Matteo Caggio}
\date{Department of Mathematics, Faculty of Science, University of Zagreb \\ 
\small matteo.caggio@math.hr}

\begin{document}

\maketitle

\begin{abstract}
    In this paper we study the incompressible inviscid limit for a compressible micro-polar model. We prove that the weak solution of the compressible micro-polar system converges to the solution of the Navier-Stokes equations (Euler equations) in the limit of small Mach number (and vanishing viscosity).
\end{abstract}


\tableofcontents{}

\newpage{}

\section{Introduction}
We consider the compressible Navier-Stokes system describing the motion of a barotropic micro-polar fluid in the whole space $\mathbb{R}^3$
\begin{equation} \label{cont}
    \partial_t \varrho_\varepsilon + \text{div}(\varrho_\varepsilon \mathbf{u}_\varepsilon)=0,
\end{equation}

\begin{equation} \label{mom}
    \partial_t (\varrho_\varepsilon \mathbf{u}_\varepsilon) + \text{div}(\varrho_\varepsilon \mathbf{u}_\varepsilon \otimes \mathbf{u}_\varepsilon) + \nabla p(\varrho_\varepsilon)
    = \text{div}\mathbb{S}(\nabla \mathbf{u}_\varepsilon) + 2 (\xi \nabla \times \boldsymbol{\omega}_\varepsilon),
\end{equation}

\begin{equation} \label{omega}
    \mathbb{I}\left[ \partial_t (\varrho_\varepsilon \boldsymbol{\omega}_\varepsilon) + \text{div}(\varrho \mathbf{u}_\varepsilon \otimes \boldsymbol{\omega}_\varepsilon)\right] = \text{div}\mathbb{M}(\nabla \boldsymbol{\omega}_\varepsilon) + 2(\xi \nabla \times \mathbf{u}_\varepsilon - 2\xi \boldsymbol{\omega}_\varepsilon).
\end{equation}
Here, $\varrho_\varepsilon=\varrho_\varepsilon(x,t)$, $\mathbf{u}_\varepsilon=\mathbf{u}_\varepsilon(x,t)$, $\boldsymbol{\omega}_\varepsilon=\boldsymbol{\omega}_\varepsilon(x,t)$ and $p=p(\varrho_\varepsilon)$ denote the density, velocity, micro-rotational (angular) velocity and pressure of the fluid, respectively. The quantity $\mathbb{I}$ represents the micro-inertia coefficient and the viscous stress tensors read as  
\begin{equation} \label{S}
    \mathbb{S}(\nabla \mathbf{u}_\varepsilon) = 
    (\mu + \xi) \nabla \mathbf{u}_\varepsilon + (\mu + \lambda - \xi) \text{div} \mathbf{u}_\varepsilon,
\end{equation}
\begin{equation} \label{M}
    \mathbb{M}(\nabla \boldsymbol{\omega}_\varepsilon) = 
    \mu' \nabla \boldsymbol{\omega}_\varepsilon + (\mu'+ \lambda') \text{div} \boldsymbol{\omega}_\varepsilon,
\end{equation}
with coefficients of viscosity $\mu, \lambda$ and of micro-viscosity $\mu', \lambda', \xi$ satisfying
\begin{equation} \label{viscosity}
    \mu, \mu', \xi >0, \ \ 2\mu + 3\lambda \geq 0, \ \ 2\mu' + 3\lambda' \geq 0.
\end{equation}

Since the works in \cite{Er} and \cite{Lu}, the compressible micro-polar system has recently received a considerable attention in terms of mathematical analysis.
In one-dimension, the compressible, viscous and heat-conducting micropolar system
was analyzed in \cite{Mu}.
The global existence of weak solutions in three-dimension with discontinuous initial data was studied in \cite{ChXuZh} and, more recently, the existence of dissipative measure-valued solutions was proved in \cite{Hu} together with a weak-strong uniqueness principle. Existence and large time behavior of strong solutions in three-dimension was concerned by \cite{LiZh}.
For further results on weak and strong solutions to the micropolar system we also refer to \cite{Am}, \cite{ChHuZh}, \cite{ChMi} \cite{Dr}, \cite{DrMu0},  \cite{Du}, \cite{FeZh}, \cite{HuLiZh}, \cite{Wu}, \cite{Ya} and references therein.
In particular, in the context of singular limit analysis, recent results are given in \cite{Su}, \cite{Su-1}, \cite{Su-2}.

For the purpose of our analysis, we rewrite the above system in non-dimensional form, namely for each physical quantity $X$ present in the Navier-Stokes equations, we introduce its reference value $X_{ref}$ and replace $X$ with its dimensionless analogue $X/X_{ref}$. 
Consequently, the system (\ref{cont}) - (\ref{omega}) reads as follows
\begin{equation} \label{cont-scl}
    \mathcal{S}r \partial_t \varrho_\varepsilon + \text{div}(\varrho_\varepsilon \mathbf{u}_\varepsilon)=0,
\end{equation}

\begin{equation} \label{mom-scl}
    \mathcal{S}r \partial_t (\varrho_\varepsilon \mathbf{u}_\varepsilon) + \text{div}(\varrho_\varepsilon \mathbf{u}_\varepsilon \otimes \mathbf{u}_\varepsilon) + \frac{\nabla p(\varrho_\varepsilon)}{\mathcal{M}a^2}
    = \frac{1}{\mathcal{R}e}\text{div}\mathbb{S}(\nabla \mathbf{u}_\varepsilon) + 2 \frac{1}{\mathcal{R}e}(\xi \nabla \times \boldsymbol{\omega}_\varepsilon),
\end{equation}

\begin{equation} \label{omega-scl}
    \partial_t (\varrho_\varepsilon \boldsymbol{\omega}_\varepsilon) + \text{div}(\varrho \mathbf{u}_\varepsilon \otimes \boldsymbol{\omega}_\varepsilon) = \frac{\mathcal{R}e_M}{\mathcal{R}e}\text{div}\mathbb{M}(\nabla \boldsymbol{\omega}_\varepsilon) + 2\frac{1}{\mathcal{R}e}(\xi \nabla \times \mathbf{u}_\varepsilon - 2\xi \boldsymbol{\omega}_\varepsilon).
\end{equation}
Here, $\mathcal{M}a$ is the Mach
number and $\mathcal{R}e$ is the Reynolds number defined as follows
$$
\mathcal{M}a = u_{ref} / \sqrt{p_{ref}/\varrho_{ref}}, \ \ \mathcal{R}e = \varrho_{ref} u_{ref} L_{ref} / \mu_{ref},
$$
while
$$
\mathcal{S}r = \frac{L_{ref}}{t_{ref} u_{ref}}
$$
is the Strouhal number. We also introduce
$$
\mathcal{R}e_M = \mu'_{ref} / \mu_{ref} L^2_{ref}.
$$
The Mach number is the ratio of the characteristic velocity of the flow to the speed of the sound in the fluid while the Reynolds number is the ratio of the inertial to the viscous forces in the fluid. The Strouhal number plays a role in oscillating, non-steady flows, as the K\' arm\' an vortex street (see e.g. \cite{Wa}), where $1/t_{ref}$ is the frequency of vortex shedding in the wake 
of von K\' arma\' n (the characteristic length $L_{ref}$ could refers to a body invested by the flow with characteristic velocity $u_{ref}$).
Here and hereafter, the Strouhal number is set equal to one and we keep 
$\mathcal{R}e_M$ fixed and positive.
For more details about the scaling procedure, the reader is referred to \cite{FeNo}, Chapter 4.

\begin{rem} \label{scal}
Note that in performing the scaling analysis we considered $\mu_{ref}$ as the reference viscosity for $\mu$, $\lambda$ and $\xi$, and $\mu'_{ref}$ for $\mu'$ and $\lambda'$. Here, $\omega_{ref} = 1/t_{ref}$ and $\mathbb{I}_{ref} = L^2_{ref}$. 
For the definition of $\mathcal{R}e_M$ see e.g. \cite{Ba} and the references therein.
\end{rem}

For simplicity of the notation, we set $\mathcal{M}a = \varepsilon$ and $\mathcal{R}e = \nu^{-1}$. 
We assume that the pressure is given by the following relation
\begin{equation} \label{pressure}
    p(\varrho_\varepsilon)= a \varrho_\varepsilon^\gamma, \ \ a >0, \ \ \gamma > 3/2.
\end{equation}
The system is supplemented with the initial condition
\begin{equation} \label{ic}
    \varrho_\varepsilon(0, \cdot)=\varrho_{0,\varepsilon} = 1 + \varepsilon \varrho_{0,\varepsilon}^{(1)}, \ \ \mathbf{u}_\varepsilon(0, \cdot)=\mathbf{u}_{0,\varepsilon}, \ \ \boldsymbol{\omega}_\varepsilon(0, \cdot)=\boldsymbol{\omega}_{0,\varepsilon} 
\end{equation}
and the far field conditions
\begin{equation} \label{bc}
    \varrho_\varepsilon \rightarrow 1, \ \ \mathbf{u}_\varepsilon \rightarrow 0, \ \ \boldsymbol{\omega}_\varepsilon \rightarrow 0 \ \ \text{as} \ \ |x| \rightarrow \infty.
\end{equation}

\begin{rem}\label{ic-r}
The initial condition (\ref{ic})${_1}$ can be seen as a perturbation of the density field respect to the constant (equilibrium) density as first order expansion in terms of the Mach number (see \cite{FeNo}, Chapter 4).
The role of the quantity $\varrho_{0,\varepsilon}^{(1)}$ will be clear later when we introduce the acoustic system related to the compressible Navier-Stokes equations (see Section 3).
\end{rem}

Our aim is to identify the system of equations in the limit of $\varepsilon\rightarrow0$ and in the limit $\varepsilon,\nu\rightarrow0$, respectively.
More precisely, in the first case our goal is to prove the convergence of the weak solution of the compressible micro-polar system to the weak solution of the incompressible Navier-Stokes equation for micro-polar fluids, namely
\begin{equation} \label{cont-NS}
    \text{div} \mathbf{v} =0,
\end{equation}
\begin{equation} \label{mom-NS}
    \partial_t \mathbf{v} +\mathbf{v}\cdot\nabla\mathbf{v} + \nabla \Pi 
    - \frac{1}{\mathcal{R}e}\text{div}\mathbb{S}(\nabla \mathbf{v})
    -2\frac{1}{\mathcal{R}e}\xi \nabla \times \boldsymbol{\omega}=0,
\end{equation}
\begin{equation} \label{micr-NS}
    \partial_t \boldsymbol{\omega} +\mathbf{v}\cdot\nabla\boldsymbol{\omega}
    - \frac{\mathcal{R}e_M}{\mathcal{R}e}\text{div}\mathbb{S}(\nabla \boldsymbol{\omega})
    +2 \frac{1}{\mathcal{R}e}(2\xi\boldsymbol{\omega}
    -\xi \nabla \times \mathbf{v})
    =0.
\end{equation}
The idea is to obtain suitable a priori estimates in order to pass to the limit in the weak formulation for $\varepsilon \rightarrow 0$. The main difficulty is to handle the convective term in the momentum equations. Namely to show that $\varrho_\varepsilon \mathbf{u}_\varepsilon \otimes \mathbf{u}_\varepsilon \rightarrow \mathbf{v} \otimes \mathbf{v}$ in a weak limit. This could be obtained through the analysis of the acoustic equations. 
The convergence is meant in the sense of the so-called \textit{ill-prepared} initial data, namely data that allow the propagation of acoustic waves (see \cite{FeNo}). The energy and Strichartz estimates (dispersive estimates) discussed below allow the dispersion of the acoustic waves in the limit $\varepsilon \rightarrow 0$. The analysis is inspired by the recent work \cite{Ca}.

In the second case we aim to prove the convergence of the weak solution of the compressible micro-polar system to the strong solution of the incompressible Euler equation for micro-polar fluids, namely 
\begin{equation} \label{cont-E}
    \text{div} \mathbf{v} =0,
\end{equation}
\begin{equation} \label{mom-E}
    \partial_t \mathbf{v} +\mathbf{v}\cdot\nabla\mathbf{v} + \nabla \Pi =0,
\end{equation}
\begin{equation} \label{micr-E}
    \partial_t \boldsymbol{\omega} +\mathbf{v}\cdot\nabla\boldsymbol{\omega}
    =0.
\end{equation}
The convergence is obtained from a Gronwall type argument introducing a \textit{relative energy inequality} (see e.g. \cite{Da}, \cite{FeNo}, \cite{Ge}, \cite{SR}). Acoustic waves are dispersed through the energy and Strichartz estimates discussed below. 

In the last years, the study of the low Mach number and/or inviscid regime for compressible fluids has been of large interest for many authors (see \cite{GiNo} and references therein). However, in the context of micro-polar compressible fluids very few results are present (see \cite{Su}, \cite{Su-1}, \cite{Su-2}). Indeed, as far as the author is aware, the weak-weak convergence in the low Mach number limit for the compressible (barotropic) micro-polar fluid in the whole space with ill-prepared initial data has not been analyzed before. Concerning the weak-strong convergence, the limit was already studied by \cite{Su}. However, the author in \cite{Su} obtains the convergence under some conditions satisfied in two ranges of the adiabatic coefficient $\gamma$, namely $3/2 < \gamma < 4$ and $\gamma \geq 4$ (see Theorem 2.3 in \cite{Su}). In the present result, we remove these conditions and generalize the result for $\gamma > 3/2$. 
We combine the recent analysis done in \cite{Hu}, in which the author proves the existence of dissipative measure-valued solution to a compressible micro-polar system together with a weak-strong uniqueness principle,  with the analysis in \cite{Ca} and \cite{CaNe} in order to treat, on one side, the dispersion of the acoustic waves, and, on the other side, the micro-polar coupling.


\subsection{Weak solution to the compressible system}

\begin{defn} \label{weak-def}
We say that $(\varrho_\varepsilon, \mathbf{u}_\varepsilon, \boldsymbol{\omega}_\varepsilon)$ is a weak solution to the compressible micro-polar Navier-Stokes system (\ref{cont}) - (\ref{bc}) if

\bigskip
$\bullet$ the functions $(\varrho_\varepsilon, \mathbf{u}_\varepsilon, \boldsymbol{\omega}_\varepsilon)$ belongs to the class
\begin{equation} \label{int-rho}
    \varrho_\varepsilon - 1 \in L^\infty (0,T; L^2(\mathbb{R}^3) + L^\gamma(\mathbb{R}^3)), \ \ \varrho_\varepsilon \geq 0 \ \ \mbox{a.a} \ \ \mbox{in} \ \ (0,T)\times\mathbb{R}^3,
\end{equation}
\begin{equation} \label{int-u}
    \mathbf{u}_\varepsilon \in L^2(0,T; D_0^{1,2}(\mathbb{R}^3;\mathbb{R}^3)), \ \ \varrho_\varepsilon \mathbf{u}_\varepsilon \in L^\infty (0,T;(L^2+L^{\frac{2\gamma}{\gamma+1}})(\mathbb{R}^3;\mathbb{R}^3)),
\end{equation}
\begin{equation} \label{int-omega}
    \boldsymbol{\omega}_\varepsilon \in L^2(0,T; D_0^{1,2}(\mathbb{R}^3;\mathbb{R}^3)), \ \ \varrho_\varepsilon \boldsymbol{\omega}_\varepsilon \in L^\infty (0,T;(L^2+L^{\frac{2\gamma}{\gamma+1}})(\mathbb{R}^3;\mathbb{R}^3)).
\end{equation}

\bigskip
$\bullet$ $\varrho_\varepsilon - 1 \in C_{weak} ([0,T];L^2(\mathbb{R}^3) + L^\gamma(\mathbb{R}^3))$ and the integral identity
\begin{equation} \label{weak-cont}
    \int_{\mathbb{R}^3} \varrho_\varepsilon(\tau, x) \varphi(\tau, x) dx - \int_{\mathbb{R}^3} \varrho_{0,\varepsilon} \varphi(0, x) dx =  \int_0^T \int_{\mathbb{R}^3} \left( \varrho_\varepsilon \partial_t \varphi + \varrho_\varepsilon \mathbf{u}_\varepsilon \cdot \nabla \varphi \right) dxdt
\end{equation}
holds for all $\tau \in [0,T]$ and any test function $\varphi \in C_c^\infty ([0,T]\times {\mathbb{R}^3}$.

\bigskip
$\bullet$ $\varrho_\varepsilon \mathbf{u}_\varepsilon \in C_{weak} \left([0,T]; (L^2 + L^{\frac{2\gamma}{\gamma+1}})(\mathbb{R}^3;\mathbb{R}^3)\right)$ and the integral identity
$$
\int_{\mathbb{R}^3} \varrho_\varepsilon \mathbf{u}_\varepsilon(\tau, x) \cdot \varphi(\tau, x) dx - \int_{\mathbb{R}^3} \varrho_{0,\varepsilon} \mathbf{u}_{0,\varepsilon} \cdot \varphi(0, x) dx
$$   
$$
= \int_0^T \int_{\mathbb{R}^3} \left( \varrho_\varepsilon \mathbf{u}_\varepsilon \cdot \partial_t \varphi + \varrho_\varepsilon \mathbf{u}_\varepsilon \otimes \mathbf{u}_\varepsilon : \nabla \varphi + \frac{p(\varrho_\varepsilon)}{\varepsilon^2} \text{div} \varphi \right) dxdt
$$
$$
- \nu \int_0^T \int_{\mathbb{R}^3} \left( (\mu + \xi) \nabla \mathbf{u}_\varepsilon : \nabla \varphi + (\mu + \lambda - \xi) \text{div} \mathbf{u}_\varepsilon \text{div} \varphi \right)dxdt
$$
\begin{equation} \label{weak-mom}
+ \nu \int_0^T \int_{\mathbb{R}^3} 2\xi \ (\nabla \times \boldsymbol{\omega}_\varepsilon) \cdot \varphi dxdt 
\end{equation}
holds for all $\tau \in [0,T]$ and any test function $\varphi \in C_c^\infty ([0,T]\times \mathbb{R}^3;\mathbb{R}^3)$.

\bigskip
$\bullet$ $\varrho_\varepsilon \boldsymbol{\omega}_\varepsilon \in C_{weak} \left([0,T]; (L^2 + L^{\frac{2\gamma}{\gamma+1}})(\mathbb{R}^3;\mathbb{R}^3)\right)$ and the integral identity
$$
\int_{\mathbb{R}^3} \varrho_\varepsilon \boldsymbol{\omega}_\varepsilon (\tau, x) \cdot \varphi(\tau, x) dx - \int_{\mathbb{R}^3} \varrho_{0,\varepsilon} \boldsymbol{\omega}_{0,\varepsilon} \cdot \varphi(0, x) dx
$$   
$$
= \int_0^T \int_{\mathbb{R}^3} \left( \varrho_\varepsilon \boldsymbol{\omega}_\varepsilon \cdot \partial_t \varphi + \varrho_\varepsilon \boldsymbol{\omega}_\varepsilon \otimes \mathbf{u}_\varepsilon : \nabla \varphi \right) dxdt
$$
$$
- \nu \int_0^T \int_{\mathbb{R}^3} \left( \mu' \nabla \boldsymbol{\omega}_\varepsilon : \nabla \varphi + (\mu' + \lambda') \text{div} \boldsymbol{\omega}_\varepsilon \text{div} \varphi \right)dxdt
$$
\begin{equation} \label{weak-omega}
+ \nu \int_0^T \int_{\mathbb{R}^3} 2\xi (\nabla \times \mathbf{u}_\varepsilon -2 \boldsymbol{\omega}_\varepsilon) \cdot \varphi dxdt 
\end{equation}
holds for all $\tau \in [0,T]$ and any test function $\varphi \in C_c^\infty ([0,T]\times \mathbb{R}^3;\mathbb{R}^3)$.

\bigskip
$\bullet$ The energy inequality
$$
\int_{\mathbb{R}^3} \left( \frac{1}{2} \varrho_\varepsilon |\mathbf{u}_\varepsilon|^2
+ \frac{1}{2} \varrho_\varepsilon |\boldsymbol{\omega}_\varepsilon|^2
+ \frac{P(\varrho_\varepsilon,1)}{\varepsilon^2}\right) (\tau, x) dx
$$
$$
+ \nu \int_0^T \int_{\mathbb{R}^3} \left( \mu|\nabla \mathbf{u}_\varepsilon|^2 +(\mu+\lambda)|\text{div}\mathbf{u}_\varepsilon|^2+\mu'|\nabla \boldsymbol{\omega}_\varepsilon|^2+(\mu'+\lambda')|\text{div}\boldsymbol{\omega}_\varepsilon|^2+\xi|2\boldsymbol{\omega}_\varepsilon-\nabla \times \mathbf{u}_\varepsilon|^2 \right) dxdt
$$
\begin{equation} \label{ei}
\leq \int_{\mathbb{R}^3} \left( \frac{1}{2} \varrho_{0,\varepsilon} |\mathbf{u}_{0,\varepsilon}|^2
+ \frac{1}{2} \varrho_{0,\varepsilon} |\boldsymbol{\omega}_{0,\varepsilon}|^2
+ \frac{P(\varrho_{0,\varepsilon},1)}{\varepsilon^2}\right) dx
\end{equation}
holds for a.a $\tau \in [0,T]$, where
$$
P(\varrho_\varepsilon,1) = H(\varrho_\varepsilon) - H'(1)(\varrho_\varepsilon-1)-H(1),
$$
with
$$
H(\varrho_\varepsilon) = \varrho_\varepsilon \int_1^{\varrho_\varepsilon} \frac{p(z)}{z^2} dz.
$$
\end{defn}

\begin{rem} \label{D_0}
Here, the space $D_{0}^{1,2}(\mathbb{R}^3)$ is a completion of $\mathcal{D}(\mathbb{R}^3)$
-- the space of smooth functions compactly supported in $\mathbb{R}^3$ --
with respect to the norm
$$
\left\Vert f\right\Vert _{D_{0}^{1,2}(\mathbb{R}^3)}^{2}=\int_{\mathbb{R}^3}\left|\nabla f\right|^{2} dx.
$$
In accordance with the Sobolev's inequality, we have
\begin{equation} \label{S-ineq}
    D_{0}^{1,2}(\mathbb{R}^3) \subset L^6(\mathbb{R}^3), 
\end{equation}
see \cite{Ga}.
\end{rem}

The following result concerns the existence of a global weak solution to the compressible micro-polar system (\ref{cont}) - (\ref{omega}) (see \cite{Su}, Proposition 2.2; see also \cite{Am}).

\begin{prop} \label{NS-weak}
    Let $\gamma > 3/2$ and the initial data $(\varrho_{0,\varepsilon},\mathbf{u}_{0,\varepsilon},\boldsymbol{\omega}_{0,\varepsilon})$ satisfy
    \begin{equation} \label{e_0}
    \int_{\mathbb{R}^3} \left( \frac{1}{2} \varrho_{0,\varepsilon} |\mathbf{u}_{0,\varepsilon}|^2
    + \frac{1}{2} \varrho_{0,\varepsilon} |\boldsymbol{\omega}_{0,\varepsilon}|^2
    + \frac{P(\varrho_{0,\varepsilon},1)}{\varepsilon^2}\right) \leq C 
    \end{equation}
    such that
    \begin{equation} \label{m_0}
        \varrho_{0,\varepsilon} \mathbf{u}_{0,\varepsilon} (x) =0, \ \ \varrho_{0,\varepsilon} \boldsymbol{\omega}_{0,\varepsilon} (x) =0, \ \ \mbox{whenever} \ \ x \in \{\varrho_{0,\varepsilon}=0\}.
    \end{equation}
    Then, there exists at least one weak solution to the Navier-Stokes system (\ref{cont}) - (\ref{omega}) in the sense of Definition \ref{weak-def}.
\end{prop}

\subsection{Weak and strong solution to the incompressible system}

We recall the global existence of weak solution for the Navier-Stokes system (\ref{cont-NS}) - (\ref{micr-NS}), and the local existence of strong solution to the incompressible inviscid system (\ref{cont-E}) - (\ref{micr-E}).

\begin{thm} \label{NS-micr}
[\cite{Lu} Chapter 3, Theorem 1.1.6]
    Assume the initial data $\mathbf{v}_0 \in W^{1,2}(\mathbb{R}^3;\mathbb{R}^3)$, $\text{div}\mathbf{v}_0=0$, and $\boldsymbol{\omega}_0 \in L^2(\mathbb{R}^3;\mathbb{R}^3)$. Then, there exist $(\mathbf{v}, \boldsymbol{\omega}, p)$
    $$
    \mathbf{v} \in L^\infty(0,T;H(\mathbb{R}^3;\mathbb{R}^3)) \cap L^2(0,T;V(\mathbb{R}^3;\mathbb{R}^3)),
    $$
    $$
    \boldsymbol{\omega} \in L^\infty(0,T;L^2(\mathbb{R}^3;\mathbb{R}^3)) \cap L^2(0,T;D_0^1(\mathbb{R}^3;\mathbb{R}^3)),
    $$
    \begin{equation} \label{reg-NS}
    p \in \mathcal{D}' (\mathbb{R}^3 \times (0,T))
    \end{equation}
    satisfying (\ref{cont-NS}) - (\ref{micr-NS}) in the sense of distribution and initial data weakly in $L^2(\mathbb{R}^3)$.  
\end{thm}
\begin{rem} \label{spaces}
Above, $(H, V)$ denote the closure of $\{ f \in C_0^\infty(\mathbb{R}^3;\mathbb{R}^3); \text{div}f=0 \}$ in $L^2(\mathbb{R}^3)$ and $D_0^1(\mathbb{R}^3)$.
\end{rem}
\begin{thm} \label{E-micr}
[\cite{Su}, Proposition 2.1]
    Let s > 5/2. Assume the initial data $\mathbf{v}_0, \boldsymbol{\omega}_0 \in W^{s,2}(\mathbb{R}^3;\mathbb{R}^3)$ and $\text{div}\mathbf{v}_0=0$. Then, there exists $T^* \in (0,\infty)$ and a unique solution $(\mathbf{v},\boldsymbol{\omega}) \in L^\infty(0,T^*;W^{s,2}(\mathbb{R}^3;\mathbb{R}^3))$ to the system (\ref{cont-E}) - (\ref{micr-E}) such that
    \begin{equation} \label{reg-E}
        \sup_{0\leq t\leq T} \{\Vert(\mathbf{v},\boldsymbol{\omega})(t,\cdot)\Vert_{W^{s,2}(\mathbb{R}^3;\mathbb{R}^3)}+\Vert(\partial_t \mathbf{v},\partial_t \boldsymbol{\omega})(t,\cdot)\Vert_{W^{s-1,2}(\mathbb{R}^3;\mathbb{R}^3)}\}\leq C_T
    \end{equation}
    for any $0<T<T^*$.
\end{thm}

\begin{rem} \label{str-s}
Because the system (\ref{cont-E}) - (\ref{micr-E}) is decoupled, the classical results on the incompressible Euler
equations (see e.g. \cite{Ka}) and on the transport equation gives the estimate (\ref{reg-E}) (see \cite{Su}).
\end{rem}

\subsection{Main results}
Our main results can be stated as follows:

\begin{thm}\label{mainweak}
Let $(\varrho_\varepsilon, \mathbf{u}_\varepsilon,\boldsymbol{\omega}_\varepsilon)$ be the weak solution to the compressible Navier-Stokes system (\ref{cont})-(\ref{omega}) with the initial data
\begin{equation}\label{bound-ic}
(\mathbf{u}_{0,\varepsilon}, \boldsymbol{\omega}_{0,\varepsilon}) \text{ bounded in }L^2(\mathbb{R}^3;\mathbb{R}^3), \ \  \varrho^{(1)}_{0,\varepsilon} \text { bounded in } (L^2 \cap L^\infty)(\mathbb{R}^3)
\end{equation}
uniformly for  $\varepsilon$ such that
\begin{equation}\label{limit-ic}
\varrho_{0,\varepsilon}\mathbf{u}_{0,\varepsilon} \to \mathbf{v}_{0},
\ \ 
\varrho_{0,\varepsilon}\boldsymbol{\omega}_{0,\varepsilon} \to \boldsymbol{\omega}_{0} \text{ in } L^{2}(\mathbb{R}^3;\mathbb{R}^3)
\end{equation}
as $\varepsilon\to 0$. Then
\begin{equation}\label{mainweak1}
\varrho_\varepsilon \to 1 \text{ in }L^\infty(0,T;(L^{2}+L^{\gamma})(\mathbb{R}^3), \ \ \mathbf{u}_{\varepsilon} \to \mathbf{v} \text{ weakly in }L^2(0,T;W_0^{1,2}(\mathbb{R}^3;\mathbb{R}^3)),
\end{equation}
\begin{equation}\label{mainweak2}
\boldsymbol{\omega}_{\varepsilon} \to \boldsymbol{\omega} \text{ weakly in }L^2(0,T;W_0^{1,2}(\mathbb{R}^3;\mathbb{R}^3))
\end{equation}
and
\begin{equation}
    \mathbf{u}_{\varepsilon} \to \mathbf{v} \text{ in }L^2(0,T;L^{2}_{loc}(\mathbb{R}^3;\mathbb{R}^3)),
\end{equation}
for any $T>0$, where $(\mathbf{v}, \boldsymbol{\omega})$ is the weak solution to the initial value problem (\ref{cont-NS})-(\ref{micr-NS}).
\end{thm}

\begin{thm}\label{maineuler}
Assume
there exist $\varrho^{(1)}_0\in L^2(\mathbb{R}^3)$
and
$
(\mathbf{u}_0, \boldsymbol{\omega}_0) \in L^2(\mathbb{R}^3;\mathbb{R}^3)$ such that
\begin{equation}\label{comp}
\left(\varrho^{(1)}_{0,\varepsilon}-{\varrho}^{(1)}_{0}\right), \ \ 
\left(\mathbf{u}_{0,\varepsilon} - \mathbf{u}_0\right), \ \ 
\left(\boldsymbol{\omega}_{0,\varepsilon} - \boldsymbol{\omega}_0\right)\rightarrow 0\text{ in } L^2(\mathbb{R}^3)
\end{equation}
and 
$(\mathbf{v}_0=\mathbf{H}(\mathbf{u}_{0}),\boldsymbol{\omega}_0)\in W^{s,2}(\mathbb{R}^3;\mathbb{R}^3)$ (s>5/2), such that $\text{supp}[\mathbf{v}_0]$ and $\text{supp}[\boldsymbol{\omega}_0]$ compact in $\mathbb{R}^3$,
$\nabla \Psi_0=\mathbf{H}^{\perp}(\mathbf{u}_{0})\in L^{2}(\mathbb{R}^3;\mathbb{R}^3)$. Let $(\mathbf{v},\boldsymbol{\omega})$ be the unique solution to the initial value problem (\ref{cont-E})-(\ref{micr-E}) and  $(\varrho_\varepsilon, \mathbf{u}_\varepsilon, \boldsymbol{\omega}_\varepsilon)$ be the weak solution to the compressible Navier-Stokes system (\ref{cont})-(\ref{omega}). Then, as $\varepsilon,\nu \rightarrow0$,
\begin{equation}\label{maineuler1}
\varrho_\varepsilon \to 1 \text{ in }L^\infty(0,T;(L^{2}+L^{\gamma})(\mathbb{R}^3)), \ \  \sqrt{\varrho_{\varepsilon}}\mathbf{u}_{\varepsilon} \to \mathbf{v} \text{ in }L^2(0,T;L^{2}_{loc}(K;\mathbb{R}^3)),
\end{equation}
\begin{equation}\label{maineuler2}
\sqrt{\varrho_{\varepsilon}}\boldsymbol{\omega}_{\varepsilon} \to \boldsymbol{\omega} \text{ in }L^2(0,T;L^{2}_{loc}(K;\mathbb{R}^3)) 
\end{equation}
for any $T>0$ and any compact set $K\subset \mathbb{R}^3$.
\end{thm}

\begin{rem} \label{Helm}
Here, the symbol $\mathbf{H}$ and $\mathbf{H}^\perp$ denote the Helmoholtz projection onto the space of solenoidal functions and the corresponding gradient part, in the sense that a vector function $\mathbf{f}: \Omega \to \mathbb{R}^3$ is written as
\begin{equation} \label{H-proj}
    \mathbf{f} = \underbrace{\mathbf{H}[\mathbf{f}]}_\text{solenoidal part}
    + 
    \underbrace{\mathbf{H}^\perp[\mathbf{f}]}_\text{gradient part}.
\end{equation}
For more details we refer the reader to \cite{FeNo}, Section 5.4.1.
\end{rem}
\begin{cor} \label{cor}
Indeed, the exact convergence rate holds
$$
\left\Vert \sqrt{\varrho_\varepsilon}\left(\mathbf{u}_\varepsilon-\mathbf{v}-\nabla\Psi_\varepsilon\right)(\tau,\cdot)\right\Vert _{L^{2}(\mathbb{R}^3;\mathbb{R}^{3})}^{2}
+
\left\Vert \sqrt{\varrho_\varepsilon}\left(\boldsymbol{\omega}_\varepsilon-\boldsymbol{\omega}\right)(\tau,\cdot)\right\Vert _{L^{2}(\mathbb{R}^3;\mathbb{R}^{3})}^{2}
$$
$$
+\left\Vert \frac{\varrho-1}{\varepsilon}(\tau,\cdot)-\psi(\tau,\cdot)\right\Vert _{L^{2}(\mathbb{R}^3)}^{2}+\left\Vert \frac{\varrho-1}{\varepsilon^{2/\gamma}}(\tau,\cdot)-\frac{\psi(\tau,\cdot)}{\varepsilon^{\left(2/\gamma\right)-1}}\right\Vert _{L^{\gamma}(\mathbb{R}^3)}^{\gamma}
$$ 
\begin{equation} \label{conv-r}
\leq c\left(\left\Vert \mathbf{u}_{0,\varepsilon}-\mathbf{u}_{0}\right\Vert _{L^{2}(\mathbb{R}^3;\mathbb{R}^{3})}^{2}+\left\Vert \varrho_{0,\varepsilon}^{(1)}-\varrho_{0}^{(1)}\right\Vert _{L^{2}(\mathbb{R}^3)}^{2}+
\left\Vert \boldsymbol{\omega}_{0,\varepsilon}-\boldsymbol{\omega}_{0}\right\Vert _{L^{2}(\mathbb{R}^3;\mathbb{R}^{3})}^{2}\right),
\end{equation}
for $\tau\in\left[0,T\right].$
\end{cor}

The paper is organized as follows. In Section 2, we introduce the uniform bounds necessary for our analysis. In Section 3, we discuss the acoustic system related to the Navier-Stokes equations.
In Section 4, we prove the convergence of the weak solution of the compressible Navier-Stokes
system to the weak solution of the incompressible system in the limit of the Mach number that tends to zero. In Section 5, through the use of the relative energy inequality, we prove the convergence of the weak solution of the compressible Navier-Stokes system to the
classical solution of the Euler equations.

\section{Preliminaries}
In the following we introduce some preliminaries estimates useful for our analysis and we discuss the a priori bounds.
\subsection{Estimates}
Let us introduce a regularizing kernel $\chi \in C_0^\infty (\mathbb{R}^3)$, $\chi \geq 0$ and $\int_{\mathbb{R}^3} \chi dx = 1$. We define $\chi_\eta = \eta^{-3}\chi(x/\eta)$. In the present analysis, we will use the following estimates (see \cite{DeGr}, relations (2.1) and (2.2)). For $f \in D^1(\mathbb{R}^3)$, we have
\begin{equation} \label{est1}
    \|f - f*\chi_\eta \|_{L^q(\mathbb{R}^3)}
    \leq
    c(q) \eta^{1-\sigma}
    \|\nabla f\|_{L^2(\mathbb{R}^3)}, \ \ q\in[2,6], \ \ \beta=3 \left(\frac{1}{2} - \frac{1}{q}\right).
\end{equation}
Moreover, for $1 < p_2 < p_1 < + \infty$, $s \geq 0$ and $\eta \in (0,1)$, we have
\begin{equation} \label{est2}
  \|f*\chi_\eta \|_{L^{p_1}(\mathbb{R}^3)}
  \leq
  c \eta^{-s - 3(1/p_2 - 1/p_1)}
  \|f\|_{W^{-s,p_2}(\mathbb{R}^3)}.
\end{equation}
\subsection{A priori bounds}
From the energy inequality, we have
\begin{equation} \label{bound_mom}
\sqrt{\varrho_{\varepsilon}}\mathbf{u}_{\varepsilon}, \ \  
\sqrt{\varrho_{\varepsilon}}\boldsymbol{\omega}_{\varepsilon}
\mbox{ uniformly bounded in } L^{\infty}\left(0,T;L^{2}\left(\mathbb{R}^3;\mathbb{R}^{3}\right)\right),
\end{equation}
\begin{equation} \label{bound_uo}
\nabla \mathbf{u}_{\varepsilon}, \ \ \nabla \boldsymbol{\omega}_{\varepsilon} \mbox{ uniformly bounded in } L^{2}\left(0,T;L^{2}\left(\mathbb{R}^3;\mathbb{R}^3\right)\right).
\end{equation}
Moreover, we observe that the map $\varrho\rightarrow P(\varrho,1)$ is
a strictly convex function on $\left(0,\infty\right)$
with global minimum equal to $0$ at $\varrho=1$ that grows at
infinity with the rate $\varrho^{\gamma}$. Consequently, the integral
$\int_{\mathbb{R}^3}P\left(\varrho,1\right)\left(\tau,\cdot\right)dx$ in (\ref{ei})
provides a control of $\left(\varrho-1\right)\left(\tau,\cdot\right)$
in $L^{2}$ over the sets $\left\{ x:\left|\varrho-1\right|\left(\tau,x\right)<\frac{1}{2}\right\} $
and in $L^{\gamma}$ over the sets $\left\{ x:\left|\varrho-1\right|\left(\tau,x\right)\geq\frac{1}{2}\right\}$.
Consequently, there holds
\begin{equation} \label{P}
P(\varrho,r)\approx\left|\varrho-1\right|^{2}1_{\left\{ \left|\varrho-r\right|<\frac{1}{2}\right\} }+\left|\varrho-r\right|^{\gamma}1_{\left\{ \left|\varrho-1\right|\geq\frac{1}{2}\right\} },\;\;\;\forall\varrho\geq0,
\end{equation}
in the sense that $P(\varrho,1)$ gives an upper and lower bound in term of the right-hand side quantity (see \cite{BaNu}, \cite{Su}). Therefore, we have the following uniform bounds

\begin{equation} \label{unif_bound1}
\textrm{ess}\,\sup_{\tau\in\left[0,T\right]}\left\Vert \left[\left(\varrho-1\right)(\tau,\cdot)\right]1_{\left\{ \left|\varrho-1\right|<\frac{1}{2}\right\} }\right\Vert _{L^{2}(\mathbb{R}^3)}\leq C\varepsilon,
\end{equation}

\begin{equation} \label{unif_bound2}
\textrm{ess}\,\sup_{\tau\in\left[0,T\right]}\left(\left\Vert \left[\left(\varrho-1\right)(\tau,\cdot)\right]1_{\left\{ \left|\varrho-1\right|\geq\frac{1}{2}\right\} }\right\Vert _{L^{\gamma}(\mathbb{R}^3)}\right)\leq C\varepsilon^{2/\gamma}.
\end{equation}
Consequently,
\begin{equation} \label{rho_conv}
\varrho_{\varepsilon}\rightarrow{1} \mbox{ in } L^{\infty}\left(0,T;L^{2}\left(\mathbb{R}^3\right)+L^{\gamma}\left(\mathbb{R}^3\right)\right).
\end{equation}
Similarly to above, for a function $f_\varepsilon$ such that $ f_\varepsilon \in D_0^{1,2}(\mathbb{R}^3)$,
we can write
\begin{equation} \label{dec-1}
    f_\varepsilon = f_\varepsilon1_{\left\{ \left|\varrho-1\right|<\frac{1}{2}\right\} } 
    + f_\varepsilon1_{\left \{ \left|\varrho-1\right|\geq\frac{1}{2} \right\}}.
\end{equation}
Now, for $f_\varepsilon = (\mathbf{u}_\varepsilon, \boldsymbol{\omega}_\varepsilon)$, following \cite{FeGaNo} we decompose
\begin{equation} \label{dec-2}
    \int_{\mathbb{R}^3} |f_\varepsilon|^2 dx = \int_{\mathbb{R}^3} (1- \varrho_\varepsilon) |f_\varepsilon|^2 dx
+ \int_{\mathbb{R}^3} \varrho_\varepsilon |f_\varepsilon|^2 dx.
\end{equation}
According to (\ref{bound_mom}), the second term in (\ref{dec-2}) is bounded in $L^\infty(0,T)$. We consider the first term and we write
\begin{equation} \label{dec-3}
    \varrho_\varepsilon - 1 = \varrho_\varepsilon^{(1)} + \varrho_\varepsilon^{(2)}
\end{equation}
with
\begin{equation} \label{conv-12}
    \varrho_\varepsilon^{(1)} \to 0 \mbox{ in } L^{\infty}\left(0,T;L^{\gamma}\left(\mathbb{R}^3\right)\right), \ \ 
    \varrho_\varepsilon^{(2)} \to 0 \mbox{ in } L^{\infty}\left(0,T;L^{2}\left(\mathbb{R}^3\right)\right),
\end{equation}
according to (\ref{unif_bound1}) and (\ref{unif_bound2}). By H\" older and Sobolev inequality, we have
$$
\left| \int_{\mathbb{R}^3} (1- \varrho_\varepsilon) |f_\varepsilon|^2 dx \right|
$$
$$
\leq c\varepsilon^{\frac{2}{\gamma}} \left(\int_{\mathbb{R}^3}\left|f_{\varepsilon}\right|^{2\gamma^{\prime}}dx\right)^{1/\gamma^{\prime}}
+ c\varepsilon \left(\int_{\mathbb{R}^3}\left|f_{\varepsilon}\right|^{4}dx\right)^{1/2}
$$
$$
\leq c\varepsilon^{\frac{2}{\gamma}}\left(\int_{\mathbb{R}^3}\left|f_{\varepsilon}\right|^{2}dx\right)^{\frac{3}{2\gamma^{\prime}}
-\frac{1}{2}}\left(\int_{\mathbb{R}^3}\left|f_{\varepsilon}\right|^{6}dx\right)^{\frac{1}{2}-\frac{1}{2\gamma^{\prime}}}
$$
$$
+ c\varepsilon \left(\int_{\mathbb{R}^3}\left|f_{\varepsilon}\right|^{2}dx\right)^{1/4}\left(
\int_{\mathbb{R}^3}\left|f_{\varepsilon}\right|^{6}dx\right)^{1/4}
$$
$$
\leq c\varepsilon^{\frac{2}{\gamma}}\left(\int_{\mathbb{R}^3}\left|f_{\varepsilon}\right|^{2}dx\right)^{\frac{3}{2\gamma^{\prime}}-\frac{1}{2}}
\left(\int_{\mathbb{R}^3}\left|\nabla f_{\varepsilon}\right|^{2}dx\right)^{\frac{3}{2}\left(1-\frac{1}{\gamma^{\prime}}\right)}
$$
\begin{equation}\label{f-essres}
+ c\varepsilon \left(\int_{\mathbb{R}^3}\left|f_{\varepsilon}\right|^{2}dx\right)^{1/4}\left(
\int_{\mathbb{R}^3}\left|\nabla f_{\varepsilon}\right|^{2}dx\right)^{3/4}
\end{equation}
with $1/\gamma + 1/\gamma' = 1$, and where we used (\ref{unif_bound1}), (\ref{unif_bound2}) and the Sobolev's inequality (\ref{S-ineq}).
By applying Young's inequality, from (\ref{bound_uo}) and (\ref{dec-2}) we conclude
\begin{equation} \label{bound_f}
f_{\varepsilon} \mbox{ uniformly bounded in } L^{2}\left(0,T;W_0^{1,2}\left(\mathbb{R}^3;\mathbb{R}^3\right)\right)
\end{equation}
and, consequently,
\begin{equation} \label{bound_grad-uo}
(\mathbf{u}_{\varepsilon},\boldsymbol{\omega}_{\varepsilon})  \mbox{ uniformly bounded in } L^{2}\left(0,T;W_0^{1,2}\left(\mathbb{R}^3;\mathbb{R}^3\right)\right),
\end{equation}
for fixed $\mu, \lambda, \mu', \lambda' >0$.
In particular, a direct consequence of (\ref{bound_mom})
gives
\begin{equation}\label{ess}
    (\mathbf{u}_\varepsilon, \boldsymbol{\omega}_\varepsilon)1_{\left\{ \left|\varrho-1\right|<\frac{1}{2}\right\} }
    \mbox{ uniformly bounded in }L^\infty\left(0,T;L^2\left(\mathbb{R}^3;\mathbb{R}^3\right)
    \right),
\end{equation}
and, back to (\ref{f-essres}), we have
\begin{equation}\label{res}
\varepsilon^{-\text{min}(1,2/\gamma)}(\mathbf{u}_{\varepsilon}, \boldsymbol{\omega}_\varepsilon)1_{\left\{ \left|\varrho-1\right|\geq\frac{1}{2}\right\} }
\mbox{ uniformly bounded in }L^2\left(0,T;L^2\left(\mathbb{R}^3;\mathbb{R}^3\right)
\right),
\end{equation}
We conclude saying that, as a consequence of (\ref{bound_mom}) and (\ref{unif_bound2}), we have
\begin{equation} \label{mom-uo-1}
    (\varrho_\varepsilon\mathbf{u}_\varepsilon, \varrho_\varepsilon \boldsymbol{\omega}_\varepsilon)1_{\left\{ \left|\varrho-1\right|<\frac{1}{2}\right\} }
    \mbox{ uniformly bounded in }L^\infty(0,T;L^2\left(\mathbb{R}^3;\mathbb{R}^3\right))
\end{equation}
and
\begin{equation} \label{mom-uo-2}
    (\varrho_\varepsilon\mathbf{u}_\varepsilon, \varrho_\varepsilon \boldsymbol{\omega}_\varepsilon)1_{\left\{ \left|\varrho-1\right|\geq\frac{1}{2}\right\} }
    \mbox{ uniformly bounded in }L^\infty(0,T;L^\frac{2\gamma}{\gamma+1}\left(\mathbb{R}^3;\mathbb{R}^3\right)).
\end{equation}

\section{Acoustic waves}

The compressibility of the fluid allows the propagation of the acoustic waves. In the following, we introduce the acoustic system related to the Eqs. (\ref{cont}) - (\ref{omega}) and we discuss the dispersion of the acoustic waves in the incompressible limit.

\subsection{Acoustic system}

The acoustic system related to the Eqs. (\ref{cont}) - (\ref{omega}) in its homogeneous form reads as follows
\begin{equation}\label{acw_1}
\varepsilon\partial_{t}{\psi}_\varepsilon + \Delta \Psi_\varepsilon=0, \ \ \varepsilon\partial_t\nabla \Psi_\varepsilon + {a}^2\nabla {\psi}_\varepsilon = 0, \ \ {a}^2=p'(1) > 0,
\end{equation}
supplemented with the initial data
\begin{equation}\label{acw_2}
{\psi}_\varepsilon(0,x) = \psi_{0, \varepsilon}(x), \ \ \nabla \Psi_\varepsilon(0,x) = \nabla \Psi_{0,\varepsilon}(x).
\end{equation}

\begin{rem} \label{ac-syst}
The homogeneous acoustic system (\ref{acw_1}) is understood as a linearization of the inviscid compressible momentum equation (\ref{mom}) where we assumed that the perturbation of the density is small in comparison to the basic state. For more details the reader can refer to \cite{Fa}, Section 2.3.1, and \cite{LaLi}, Chapter 8.
\end{rem}

\subsection{Initial data}

We consider 

\begin{equation} \label{acw}
\psi_{0,\varepsilon} = \varrho_{0,\varepsilon}^{(1)}, \ \  \nabla \Psi_{0,\varepsilon} = \mathbf{H}^{\perp}(\mathbf{u}_{{0,\varepsilon}}).
\end{equation}
The initial data (\ref{acw}) has to be understood in the sense that
the velocity $\mathbf{u}_\varepsilon$ is written in terms of the Helmholtz decomposition
$$
\mathbf{u}_{\varepsilon}=\mathbf{H}\left[\mathbf{u}_{\varepsilon}\right]+\mathbf{H}^{\perp}\left[\mathbf{u}_{\varepsilon}\right],
$$
where
$$
\mathbf{H}^{\perp}\left[\mathbf{u}_{\varepsilon}\right]=\nabla \Psi_{\varepsilon}
$$
represents the presence of the acoustic waves, with $\Psi_\varepsilon$ the acoustic potential, disappearing in the limit $\varepsilon\rightarrow0$, and $\mathbf{H}\left[\mathbf{u}_{\varepsilon}\right]$ the solenoidal component.

\subsection{Energy and decay estimates}\label{strichart}

The acoustic system conserves energy, namely
\begin{equation}\label{acwenergy}
\frac{1}{2}\int_{\mathbb{R}^3} \left|a\psi_{\varepsilon}(t,x)\right|^2 + \left|\nabla \Psi_{\varepsilon}(t,x)\right|^2 dx = \frac{1}{2}\int_{\mathbb{R}^3} \left|a\psi_0(x)\right|^2 + \left|\nabla \Psi(x)\right|^2 dx
\end{equation}
for any $t\ge 0$. Moreover, the following energy estimates hold
\[
\| \psi_\varepsilon(t,\cdot)\|_{W^{k,2}(\mathbb{R}^3)} + \|\nabla \Psi_{\varepsilon}(t,\cdot)\|_{W^{k,2}(\mathbb{R}^3;\mathbb{R}^3)}
\]
\begin{equation}\label{acw_3}
\leq c \left( \| \psi_{0,\varepsilon} \|_{W^{k,2}(\mathbb{R}^3)} + \|\nabla \Psi_{0,\varepsilon}\|_{W^{k,2}(\mathbb{R}^3;\mathbb{R}^3)} \right)
\end{equation}
for $k=1,2,\cdots,m$. Concerning the decay of the acoustic waves in the incompressible limit, we recall the following Strichartz or dispersive estimates (see \cite{DeGr}, \cite{St})
\[
\|\psi_\varepsilon\|_{L^q(0,T;L^{p}(\mathbb{R}^3))} + \|\nabla \Psi_\varepsilon\|_{L^q(0,T;L^{p}(\mathbb{R}^3;\mathbb{R}^3))}
\]
\begin{equation}\label{acw_60}
\leq c\varepsilon^{\frac{1}{q}}\left(\|\psi_{0,\varepsilon}\|_{W^{\sigma,2}(\mathbb{R}^3)}+ \|\nabla \Psi_{0,\varepsilon}\|_{W^{\sigma,2}(\mathbb{R}^3;\mathbb{R}^3)}\right),
\end{equation}
for  any
\begin{equation}\label{pqrelation}
p,q\in (2,\infty), \, \frac{1}{q}=\frac{1}{2}-\frac{1}{p},
\, \sigma = \frac{2}{q}<1.
\end{equation}
Hence for any $k=0,1,\cdots, m-1$,
\[
\|\psi_\varepsilon\|_{L^q(0,T;W^{k,p}(\mathbb{R}^3))} + \|\nabla \Psi_\varepsilon\|_{L^q(0,T;W^{k,p}(\mathbb{R}^3;\mathbb{R}^3))}
\]
\begin{equation}\label{acw_6}
\leq c\varepsilon^{\frac{1}{q}}\left(\|\psi_{0,\varepsilon}\|_{W^{m,2}(\mathbb{R}^3)}+ \|\nabla \Psi_{0,\varepsilon}\|_{W^{m,2}(\mathbb{R}^3;\mathbb{R}^3)}\right).
\end{equation}

The homogeneous form of the energy and dispersive estimates discussed above and related to the system (\ref{acw_1}) - (\ref{acw_2})
will be used in order to prove the weak-strong convergence result.

\subsubsection{Inhomogeneous case}
Now, we consider the inhomogeneous case of (\ref{acw_1}), namely
\begin{equation}\label{iacw_1}
\varepsilon\partial_{t}{\psi}_\varepsilon + \Delta \Psi_\varepsilon=\varepsilon f_1, \ \  \varepsilon\partial_t\nabla \Psi_\varepsilon + {a}^2\nabla {\psi}_\varepsilon = \varepsilon \mathbf{f}_2
\end{equation}
supplemented with the initial data
\begin{equation}\label{iacw_2}
{\psi}_\varepsilon(0,x) = \psi_{0,\varepsilon}(x), \ \  \nabla\Psi_\varepsilon(0,x) = \nabla\Psi_{0,\varepsilon}(x),
\end{equation}
where $f_1,\mathbf{f}_2\in L^q(0,T;W^{m,2}(\mathbb{R}^3))$ and $ \psi_{0,\varepsilon}(x_h),\nabla\Psi_{0,\varepsilon}\in W^{m,2}(\mathbb{R}^3)$. The following energy estimates hold
\[
\|\psi_\varepsilon\|_{L^\infty(0,T;W^{k,2}(\mathbb{R}^3))} + \|\nabla \Psi_\varepsilon\|_{L^\infty(0,T;W^{k,2}(\mathbb{R}^3;\mathbb{R}^3))}
\]
\[
\leq c\left(\|\psi_{0,\varepsilon}\|_{W^{m,2}(\mathbb{R}^3)}+ \|\nabla \Psi_{0,\varepsilon}\|_{W^{m,2}(\mathbb{R}^3;\mathbb{R}^3)}\right)
\]
\begin{equation}\label{acw_61}
+ c\left(\|f_1\|_{L^2(0,T;W^{m,2}(\mathbb{R}^3))}+ \|\mathbf{f}_2\|_{L^2(0,T;W^{m,2}(\mathbb{R}^3;\mathbb{R}^3))}\right),
\end{equation}
as well as the Strichartz or dispersive estimates
\[
\|\psi_\varepsilon\|_{L^q(0,T;W^{k,p}(\mathbb{R}^3))} + \|\nabla \Psi_\varepsilon\|_{L^q(0,T;W^{k,p}(\mathbb{R}^3;\mathbb{R}^3))}
\]
\[
\leq c\varepsilon^{\frac{1}{q}}\left(\|\psi_{0,\varepsilon}\|_{W^{m,2}(\mathbb{R}^3)}+ \|\nabla \Psi_{0,\varepsilon}\|_{W^{m,2}(\mathbb{R}^3;\mathbb{R}^3)}\right)
\]
\begin{equation}\label{iacw_4}
+ c(T)\varepsilon^{\frac{1}{q}}\left(\|f_1\|_{L^q(0,T;W^{m,2}(\mathbb{R}^3))}+ \|\mathbf{f}_2\|_{L^q(0,T;W^{m,2}(\mathbb{R}^3;\mathbb{R}^3))}\right)
\end{equation}
for the same $k,p,q$ as above (see \cite{Ca}, \cite{DeGr}).

In the weak-weak convergence analysis, we will consider the inhomogeneous case discussed above in terms of the Lighthill acoustic analogy \cite{Li-1}, \cite{Li-2}. In this sense,
the momentum $\varrho_\varepsilon \mathbf{u}_\varepsilon
=\mathbf{m}_{\varepsilon}$ will be written in terms of its Helmholtz decomposition, namely
$$
\mathbf{m}_{\varepsilon}=\mathbf{H}\left[\mathbf{m}_{\varepsilon}\right]+\mathbf{H}^{\perp}\left[\mathbf{m}_{\varepsilon}\right],
$$
where
$$
\mathbf{H}^{\perp}\left[\mathbf{m}_{\varepsilon}\right]=\nabla \Psi_{\varepsilon}.
$$
In the following we will show the compactness of the solenoidal part and that $\nabla \Psi_{\varepsilon}$ tends to zero on compact subsets and therefore becames negligible in the limit $\varepsilon \rightarrow 0$.

\subsubsection{Regularization}
For the purpose of our analysis, it is convenient to regularize the initial data (\ref{acw}) in the following way

\begin{equation} \label{smooth}
\varrho_{0, \varepsilon,\eta}^{(1)}=\chi_{\eta}\star\left(\phi_{\eta}\varrho_{0,\varepsilon}^{(1)}\right), \ \ 
\nabla\Psi_{0,\varepsilon,\eta}=\chi_{\eta}\star\left(\phi_{\eta}\nabla\Psi_{0,\varepsilon}\right),
\end{equation}
with $\eta > 0$ and where $\left\{ \chi_{\eta}\right\} $ is a family of regularazing
kernels (see Section 2.1) and $\phi_{\eta}\in C_{0}^{\infty}(\mathbb{R}^{3})$ are standards
cut-off functions. Consequently, the acoustic system posses a (unique)
smooth solution $\left[\psi_\varepsilon,\Psi_\varepsilon \right]$ and the quantities $\nabla\Psi_\varepsilon$
and $\psi_\varepsilon$ are compactly supported in $\mathbb{R}^{3}$.

\section{Weak--weak limit}

Following Lighthill \cite{Li-1}, \cite{Li-2}, it is possible to write the Navier-Stokes system in its acoustic analogy, namely
\begin{equation} \label{ac_an-1}
\varepsilon\partial_{t}\left(\frac{\varrho_{\varepsilon}-{1}}{\varepsilon}\right)+\textrm{div}\left(\varrho_{\varepsilon}\mathbf{u}_{\varepsilon}\right)=0,
\end{equation}
$$\varepsilon\partial_{t}\left(\varrho_{\varepsilon}\mathbf{u}_{\varepsilon}\right) + a^2\nabla\left(\frac{\varrho_{\varepsilon}-{1}}{\varepsilon}\right)
$$
\begin{equation*}
=\varepsilon\left( \nu\textrm{div}\mathbb{S}(\nabla\mathbf{u}_{\varepsilon})-\textrm{div}\varrho_{\varepsilon}\mathbf{u}_{\varepsilon}\otimes\mathbf{u}_{\varepsilon}\right.
\end{equation*}
\begin{equation} \label{ac_an-2}
\left.-\frac{1}{\varepsilon^{2}}\nabla \left(p\left(\varrho_{\varepsilon}\right)-a^2\left(\varrho_{\varepsilon}-{1}\right)-p\left({1}\right)\right)+2\nu\xi\nabla \times \boldsymbol{\omega}_\varepsilon \right)
\end{equation}
supplemented with the conditions (\ref{bc}).
The system (\ref{ac_an-1}) and (\ref{ac_an-2}) has to be understood in the weak sense, namely
\begin{equation} \label{ac_an-1_weak}
\int_{0}^{T}\int_{\mathbb{R}^3}\varepsilon{r_{\varepsilon}}\partial_{t}\varphi+\mathbf{m}_{\varepsilon}\cdot\nabla{\varphi}dxdt + \varepsilon \int_{\mathbb{R}^3} r_{0,\varepsilon} \varphi(0,x) dx =0,
\end{equation}
holds for every $\varphi\in C_{c}^{\infty}\left(\left[0,T\right)\times\mathbb{R}^3\right)$, while
$$
\int_{0}^{T}\int_{\mathbb{R}^3}\varepsilon\mathbf{m}_{\varepsilon}\cdot\partial_{t}\varphi+a^2{r_{\varepsilon}}\textrm{div}\varphi dxdt + \varepsilon \int_{\mathbb{R}^3}\mathbf{m}_{0,\varepsilon}\cdot\varphi(0,x) dx
$$
\begin{equation} \label{ac_an-2_weak}
=  \varepsilon\int_{0}^{T}\int_{\mathbb{R}^3} (\mathbf{f}_{\varepsilon} :\nabla\varphi + \mathbf{g}_\varepsilon\cdot \varphi) dxdt,
\end{equation}
for any $\varphi\in C_{c}^{\infty}\left(\left[0,T\right)\times\mathbb{R}^3;\mathbb{R}^{3}\right)$, where
$$
r_{\varepsilon}=\frac{\varrho_{\varepsilon}-{1}}{\varepsilon}, \ \ \mathbf{m}_{\varepsilon}=\varrho_{\varepsilon}\mathbf{u}_{\varepsilon}, \ \ \mathbf{f}_{\varepsilon} =\mathbf{f}_{\varepsilon}^{1} + \mathbf{f}_{\varepsilon}^{2} + \mathbf{f}_{\varepsilon}^{3},
$$
$$
\mathbf{f}_{\varepsilon}^{1}= \varrho_{\varepsilon}\mathbf{u}_{\varepsilon}\otimes\mathbf{u}_{\varepsilon}, \ \ \mathbf{f}_{\varepsilon}^{2}= - \nu \mathbb{S}(\nabla\mathbf{u}_{\varepsilon}),
$$
$$\mathbf{f}_{\varepsilon}^3 = \frac{1}{\varepsilon^{2}}\left(p\left({\varrho_{\varepsilon}}\right)- a^2 \left(\varrho_{\varepsilon}-{1}\right)-p\left({1}\right)\right)\mathbb{I}_2,
$$
$$
\mathbf{g}_\varepsilon =
2\nu\xi(\nabla \times \boldsymbol{\omega}_\varepsilon),
$$
such that
\begin{equation}\label{f_2}
\mathbf{f}_{\varepsilon}^{2}, \mathbf{g}_\varepsilon \mbox{ uniformly bounded in } L^{2}\left(0,T;L^{2}\left(\mathbb{R}^3;\mathbb{R}^{3\times 3}\right)\right)
\end{equation}
and $\mathbf{f}_{\varepsilon}^{1},\, \mathbf{f}_{\varepsilon}^{3} \mbox{ uniformly bounded in }  L^{\infty}(0,T;L^{1}(\mathbb{R}^3;\mathbb{R}^{3\times 3}))$ according to (\ref{bound_mom}) - (\ref{unif_bound2}).
Consequently, we have
\begin{equation} \label{f_1}
\mathbf{f}_{\varepsilon}^{1},\, \mathbf{f}_{\varepsilon}^{3} \mbox{ uniformly bounded in }  L^{\infty}(0,T;W^{-s,2}(\mathbb{R}^3;\mathbb{R}^{3\times 3})), \ \ s>3/2,
\end{equation}
since $L^{1}(\mathbb{R}^3)$ is continuously embedded in $W^{-s,2}(\mathbb{R}^3;\mathbb{R}^3)$.

\subsection{Compactness of the solenoidal component}
From the uniform bound (\ref{bound_grad-uo}), there exists $\mathbf{U}(t,x)\in\mathbb{R}^3$ such that
\begin{equation} \label{weak_conv_u}
\mathbf{u}_{\varepsilon}\rightarrow\mathbf{U} \mbox{ weakly in } L^{2}\left(0,T;W^{1,2}\left(\mathbb{R}^3;\mathbb{R}^{3}\right)\right).
\end{equation}
Consequently, from (\ref{rho_conv}), (\ref{weak_conv_u}) and the weak formulation of the continuity equation, we have
$$
\mathbf{m} = \mathbf{U}
$$
and
$$\textrm{div}\mathbf{U}=0 \mbox{ in } \mathcal{D}^{\prime},$$
which is equivalent to
$$\textrm{div}\mathbf{v}=0. $$
Moreover, from (\ref{bound_grad-uo}), we also have
\begin{equation} \label{weak_conv_omega}
\boldsymbol{\omega}_{\varepsilon}\rightarrow\boldsymbol{\omega} \mbox{ weakly in } L^{2}\left(0,T;W^{1,2}\left(\mathbb{R}^3;\mathbb{R}^{3}\right)\right).
\end{equation}
Now, we have 
\begin{equation}\label{weak_conv_m}
\mathbf{H}\left[\mathbf{m}_{\varepsilon}\right]
\rightarrow \mathbf{m} \mbox{ weakly-(*) in } L^{\infty}\left(0,T;\left(L^{2}+L^{2\gamma/\left(\gamma+1\right)}\right)\left(\mathbb{R}^3;\mathbb{R}^{3}\right)\right).
\end{equation}
From (\ref{ac_an-2_weak}) and the bounds (\ref{f_2}) and (\ref{f_1}), the above convergence could be strengthen in the sense that
\begin{equation} \label{weak_continuity}
\left[\tau\mapsto\int_{\mathbb{R}^3}\mathbf{H}\left[\mathbf{m}_{\varepsilon}\right]\cdot\phi dx\right]\rightarrow\left[\tau\mapsto\int_{\mathbb{R}^3}{\mathbf{m}}\cdot\phi dx\right] \mbox{ in } C\left[0,T\right],
\end{equation}
for any $\phi\in C_{c}^{\infty}\left(\mathbb{R}^3;\mathbb{R}^{3}\right)$, $\textrm{div}\phi=0$ (for more details see \cite{FeKaPo}, Section 6.2). 
Consequently, 
\begin{equation}\label{strong_conv_m}
\mathbf{H}\left[\mathbf{m}_{\varepsilon}\right]
\rightarrow \mathbf{m} \mbox{ in } C_{weak}\left(0,T;\left(L^{2}+L^{2\gamma/\left(\gamma+1\right)}\right)\left(\mathbb{R}^3;\mathbb{R}^{3}\right)\right).
\end{equation}
In particular, we note that $2\gamma/(\gamma+1) > 6/5$ for $\gamma > 3/2$, and the convergence (\ref{strong_conv_m}) yields
\begin{equation}\label{strong_conv_m-1}
\mathbf{H}\left[\mathbf{m}_{\varepsilon}\right]
\rightarrow \mathbf{m} \mbox{ in } C\left(0,T;W^{-1,2}\left(\mathbb{R}^3;\mathbb{R}^{3}\right)\right)
\end{equation}
thanks to the compact embedding of $L^q$ in $W^{-1,2}$ for $q > 6/5$.
Now, from (\ref{rho_conv}), (\ref{weak_conv_u}) and (\ref{strong_conv_m-1}), we have
$$
\mathbf{H}[\mathbf{u}_\varepsilon]\cdot\mathbf{H}[\mathbf{u}_\varepsilon] = \left(\mathbf{H}\left[ (1-\varrho_\varepsilon)\mathbf{u}_\varepsilon\right] + \mathbf{H}[\varrho_\varepsilon \mathbf{u}_\varepsilon]\right)\cdot \mathbf{H}[\mathbf{u}_\varepsilon] \to |\mathbf{v}|^2
$$
in the sense of distribution. 
In particular, thanks to (\ref{weak_conv_u}) and the compact embedding of $W^{1,2}(\mathbb{R}^3)$ in $L^2_{loc}(\mathbb{R}^3)$, we have
\begin{equation} \label{strong-u}
    \mathbf{H}[\mathbf{u}_{\varepsilon}] \to \mathbf{v} \text{ in }L^2(0,T;L^2_{loc}(\mathbb{R}^3;\mathbb{R}^3)).
\end{equation}

\subsection{Compactness of the gradient component}
From (\ref{ac_an-1}) - (\ref{ac_an-2}) (or from the weak formulation (\ref{ac_an-1_weak}) - (\ref{ac_an-2_weak})), we have ${r_\varepsilon} = \frac{\varrho_\varepsilon - 1}{\varepsilon}$ 
and $\nabla \Psi_{\varepsilon} = \mathbf{H}^{\perp}(\varrho_{\varepsilon}\mathbf{u}_{\varepsilon})$,
and the inhomogeneous acoustic system reads as follows
\begin{equation} \label{realac_1}
\varepsilon\partial_t {r_\varepsilon} + \Delta\Psi_\varepsilon = 0, \ \ 
\varepsilon\partial_t \nabla \Psi_\varepsilon + a^2\nabla {r_\varepsilon} = \varepsilon(\text{div}\mathbb{F}_\varepsilon+\mathbb{G}_\varepsilon),
\end{equation}
supplemented with the initial data
\begin{equation} \label{realac_2}
{r_\varepsilon}(0,x) = \varrho^{(1)}_{0,\varepsilon}, \ \ \nabla\Psi_\varepsilon(0,x)= \mathbf{H}^{\perp}(\varrho_{0,\varepsilon}\mathbf{u}_{0,\varepsilon}).
\end{equation}
Here, $\mathbb{F}_\varepsilon$ and $\mathbb{G}_\varepsilon$ are the gradient part of $\mathbf{f}_\varepsilon$ and $\mathbf{g}_\varepsilon$ in the sense of (\ref{H-proj}), namely 
$$
    \mathbb{F}_\varepsilon=
    \mathbf{H}^{\perp}
    \bigg[
    \varrho_{\varepsilon}\mathbf{u}_{\varepsilon}\otimes\mathbf{u}_{\varepsilon}
    -\nu
    \mathbb{S}(\nabla\mathbf{u}_{\varepsilon})
    +\frac{1}{\varepsilon^{2}}
    \left(p\left(\varrho_{\varepsilon}\right)-a^2\left(\varrho_{\varepsilon}-{1}\right)-p\left({1}\right)\right)\bigg]
$$
\begin{equation} \label{F-grad}
    = \mathbb{F}^1_\varepsilon
    +\mathbb{F}^2_\varepsilon
    +\mathbb{F}^3_\varepsilon,
\end{equation}

\begin{equation} \label{G-grad}
    \mathbb{G}_\varepsilon=
    \mathbf{H}^{\perp}
    \big[
    2\nu\xi\nabla \times \boldsymbol{\omega}_\varepsilon
    \big],
\end{equation}
where, in particular, we observe that $\mathbb{G}_\varepsilon=0$. Moreover,
\begin{equation}\label{FG-1}
\mathbb{F}_{\varepsilon}^{2}
\mbox{ uniformly bounded in } L^{2}\left(0,T;L^{2}\left(\mathbb{R}^3;\mathbb{R}^{3\times 3}\right)\right)
\end{equation}
\begin{equation}\label{FG-2}
\mathbb{F}_{\varepsilon}^{1},\, \mathbb{F}_{\varepsilon}^{3} \mbox{ uniformly bounded in }  L^{\infty}(0,T;W^{-s,2}(\mathbb{R}^3;\mathbb{R}^{3\times 3})), \ \ s>3/2,    
\end{equation}
according to (\ref{f_2}) and (\ref{f_1}).

The system (\ref{realac_1}) and (\ref{realac_2}) is the inhomogeneous acoustic system (\ref{iacw_1}) - (\ref{iacw_2}). In order to apply the dispersive estimates, we regularize (\ref{realac_1}) and (\ref{realac_2}) in the sense of (\ref{smooth}), namely
\begin{equation} \label{realac_1-reg}
\varepsilon\partial_t {r_{\varepsilon,\eta}} + \Delta\Psi_{\varepsilon,\eta} = 0, \ \ 
\varepsilon\partial_t \nabla \Psi_{\varepsilon,\eta} + a^2\nabla {r_{\varepsilon,\eta}} = \varepsilon
\text{div}\mathbb{F}_{\varepsilon,\eta}
\end{equation}
\begin{equation} \label{realac_2-reg}
{r_{\varepsilon,\eta}}(0,x) = \left(\varrho^{(1)}_{0,\varepsilon}\right)_\eta, \ \ \nabla\Psi_{\varepsilon,\eta}(0,x)= \left(\mathbf{H}^{\perp}(\varrho_{0,\varepsilon}\mathbf{u}_{0,\varepsilon})\right)_\eta.
\end{equation}
Now, 
we use the Strichartz estimates (\ref{iacw_4}) (with $k=0$ ($m=1$), $p=4$ and $q=4$, for example). We have
\[
\|{r_{\varepsilon}}_{,\eta}\|_{L^q_T(L^{p}(\mathbb{R}^3))} + \|\nabla\Psi_{\varepsilon,\eta}\|_{L^q_T(L^{p}(\mathbb{R}^3;\mathbb{R}^3))}
\]
\[
\leq c\varepsilon^{\frac{1}{q}}\left(\left\|{r}_{\varepsilon,\eta}(0,x)\right\|_{W^{1,2}(\mathbb{R}^3)} + \|\nabla \Psi_{\varepsilon,\eta}(0,x)\|_{W^{1,2}(\mathbb{R}^3;\mathbb{R}^3)}\right) 
\]
\begin{equation} \label{disp-reg}
+ c(T)\varepsilon^{\frac{1}{q}} \|\text{div}\mathbb{F}_{\varepsilon,\eta}\|_{L^q_T(W^{1,2}(\mathbb{R}^3;\mathbb{R}^3))}.
\end{equation}
According to the uniform bounds (\ref{FG-1}), (\ref{FG-2}) and (\ref{realac_2-reg}), the above dispersive estimate holds for 
$\text{div}\mathbb{F}_{\varepsilon}^{1}$
and $\text{div}\mathbb{F}_{\varepsilon}^{3}$. 
However, due to the lack of integrability in time, the dispersive estimate (\ref{disp-reg}) does not hold for $\text{div}\mathbb{F}_{\varepsilon}^{2}$. 
Consequently, we have
\[
\|{r_{\varepsilon}}_{,\eta}\|_{L^q_T(L^{p}(\mathbb{R}^3))} + \|\nabla\Psi_{\varepsilon,\eta}\|_{L^q_T(L^{p}(\mathbb{R}^3;\mathbb{R}^3))}
\]
\[
\leq c\varepsilon^{\frac{1}{q}}\left(\left\|{r}_{\varepsilon\\
,\eta}(0,x)\right\|_{W^{1,2}(\mathbb{R}^3)} + \|\nabla \Psi_{\varepsilon,\eta}(0,x)\|_{W^{1,2}(\mathbb{R}^3;\mathbb{R}^3)}\right) 
\]
\[
+ c(T)\varepsilon^{\frac{1}{q}} \|\text{div}\mathbb{F}^{1}_{\varepsilon,\eta}\|_{L^q_T(W^{1,2}(\mathbb{R}^3;\mathbb{R}^3))}
+ c(T)\varepsilon^{\frac{1}{q}} \|\text{div}\mathbb{F}^{2}_{\varepsilon,\eta}\|_{L^q_T(W^{1,2}(\mathbb{R}^3;\mathbb{R}^3))}
\]
\begin{equation} \label{disp-reg-1}
\leq c(\eta)\varepsilon^{\frac{1}{q}} + c(\eta,T)\varepsilon^{\frac{1}{q}}, \ \ \eta\in (0,1).
\end{equation}
Now, in order to overcome the lack of integrability in time, following the decomposition (\ref{dec-1}), in the spirit of \cite{DeGr} we define
$$
\mathbf{u}_\varepsilon^{(1)}
= 
\mathbf{u}_\varepsilon
1_{\{|\varrho-1|<\frac{1}{2}\}}, \ \ 
\mathbf{u}_\varepsilon^{(2)}
= 
\mathbf{u}_\varepsilon
1_{\{|\varrho-1|\geq\frac{1}{2}\}}. 
$$
Consequently, we split 
$$\text{div}\mathbb{F}_{\varepsilon}^{2}=\text{div}\mathbb{S}(\nabla \mathbf{u}_{\varepsilon}^{(1)})
+\text{div}\mathbb{S}(\nabla \mathbf{u}_{\varepsilon}^{(2)})
$$
$$
= \mathbb{F}_\varepsilon^{(2,1)}
+\mathbb{F}_\varepsilon^{(2,2)}.
$$
From the uniform estimates (\ref{ess}), 
we deduce that 
$$
\mathbb{F}^{(2,1)}_{\varepsilon}
\mbox{ uniformly bounded in }  L^\infty(0,T;W^{-2,2}(\mathbb{R}^3;\mathbb{R}^3)),
$$
which can be handled using (\ref{disp-reg}). While, from the uniform estimate (\ref{res}), we can deduce 
$$
\varepsilon^{-\text{min}(1,2/\gamma)}
\mathbb{F}^{(2,2)}_{\varepsilon}
\mbox{ uniformly bounded in } L^2(0,T;W^{-2,2}(\mathbb{R}^3;\mathbb{R}^3));
$$
see e.g. \cite{DeGr}.
Consequently, 
the corresponding acoustic wave produced by 
$
\mathbb{F}^{(2,2)}_{\varepsilon,\eta}
$
vanishes in $L^2(0,T;L^p(\mathbb{R}^3;\mathbb{R}^3))$ as $\varepsilon\to 0$ (for fixed $\eta$) according to (\ref{disp-reg}).
Letting $\varepsilon\to 0$, we find that for any $\eta\in (0,1)$ and $p > 2$
\begin{equation}\label{conv_grad1}
\nabla\Psi_{\varepsilon,\eta} \to 0 \text{ in }L^2(0,T;L^p
(\mathbb{R}^3;\mathbb{R}^3)),
\end{equation}
since $q>2$. 
Now, following the analysis in \cite{DeGr}, we study the following quantity
$$
\mathbf{H}^{\perp}(\mathbf{u}_{\varepsilon}) = \mathbf{H}^{\perp}(\mathbf{u}_{\varepsilon})-\mathbf{H}^{\perp}(\mathbf{u}_{\varepsilon})_\eta-\varepsilon\mathbf{H}^{\perp}\left({r}_\varepsilon\mathbf{u}_{\varepsilon}\right)_\eta + \mathbf{H}^{\perp}(\varrho_\varepsilon\mathbf{u}_{\varepsilon})_\eta
$$
where $\nabla\Psi_{\varepsilon,\eta} = \mathbf{H}^{\perp}(\varrho_\varepsilon\mathbf{u}_{\varepsilon})_\eta$.
We have
$$
\|\mathbf{H}^{\perp}(\mathbf{u}_{\varepsilon})\|_{L^2_T(L^{p}(\mathbb{R}^3;\mathbb{R}^3))} 
\leq 
\|\mathbf{H}^{\perp}(\mathbf{u}_{\varepsilon})
-\mathbf{H}^{\perp}(\mathbf{u}_{\varepsilon})_\eta\|_{L^2_T(L^{p}(\mathbb{R}^3;\mathbb{R}^3))}
$$
\begin{equation} \label{mink}
+\|\varepsilon\mathbf{H}^{\perp}\left({r}_\varepsilon\mathbf{u}_{\varepsilon}\right)_\eta\|_{L^2_T(L^{p}(\mathbb{R}^3;\mathbb{R}^3))} 
+\|\mathbf{H}^{\perp}(\varrho_\varepsilon\mathbf{u}_{\varepsilon})_\eta\|_{L^2_T(L^{p}(\mathbb{R}^3;\mathbb{R}^3))}.
\end{equation}
By the use of the estimates (\ref{est1}) and (\ref{est2}), the first and the second terms on the right hand side of (\ref{mink}) could be estimated as follows
$$
\|\mathbf{H}^{\perp}(\mathbf{u}_{\varepsilon})
-\mathbf{H}^{\perp}(\mathbf{u}_{\varepsilon})_\eta\|_{L^2_T(L^{p}(\mathbb{R}^3;\mathbb{R}^3))} \leq \eta^{1-\beta} \| \nabla \mathbf{u}_\varepsilon\|_{L^2_T(L^{2}(\mathbb{R}^3;\mathbb{R}^3))}, \ \ \beta = 3 \left(\frac{1}{2}-\frac{1}{p}\right),
$$
$$
\|\varepsilon\mathbf{H}^{\perp}\left({r}_\varepsilon\mathbf{u}_{\varepsilon}\right)_\eta\|_{L^2_T(L^{p}(\mathbb{R}^3;\mathbb{R}^3))}
\leq
\varepsilon\eta^{-(1+\beta)}\|r_\varepsilon \mathbf{u}_\varepsilon\|_{L^2_T(W^{-1,2}(\mathbb{R}^3;\mathbb{R}^3))},
$$
with $\beta \in (0,1)$ for $p \in (2,6)$. 
Consequently,
choosing $\eta$ in terms of $\varepsilon$, namely $\eta = \varepsilon^{1/\alpha}$ with $\alpha > (1+\beta)$, and taking $\varepsilon \to 0$ in (\ref{mink}), we obtain
\begin{equation}\label{conv_u_grad}
\mathbf{H}^{\perp}(\mathbf{u}_{\varepsilon}) \to 0 \text{ in }L^2(0,T;L^p(\mathbb{R}^3;\mathbb{R}^3)), \ \ p \in (2,6).
\end{equation}
Combining the convergence (\ref{conv_u_grad}) with the strong convergence of the solenoidal part (\ref{strong-u}), we obtain
\begin{equation}\label{conv_u_final}
\mathbf{u}_{\varepsilon} \to \mathbf{v} \text{ in }L^2(0,T;L^2_{loc}(\mathbb{R}^3;\mathbb{R}^3)).
\end{equation}

\subsection{The weak-weak limit convergence}
Applying the strong convergence 
(\ref{conv_u_final}) together with (\ref{weak_conv_u}) and (\ref{weak_conv_omega}) in the weak formulation, we obtain
$$
\int_{\mathbb{R}^3}\mathbf{v}\cdot\nabla\varphi dx = 0
$$
for any $\varphi\in C_{c}^{\infty}(\mathbb{R}^3)$.
Moreover,
$$
\int_{\mathbb{R}^3}\mathbf{v}\cdot\varphi(\tau,x) dx -\int_{\mathbb{R}^3} \mathbf{v}_0\cdot\varphi(0,x) dx
$$
$$
=\int_0^{\tau}\int_{\mathbb{R}^3}\mathbf{v}\cdot\partial_t\varphi + \mathbf{v}\otimes\mathbf{v}:\nabla\varphi dxdt -(\mu + \xi)\int_0^{\tau}\int_{\mathbb{R}^3}\nabla\mathbf{v}:\nabla\varphi dxdt
$$
$$
+2\xi \int_0^{\tau}\int_{\mathbb{R}^3}
(\nabla \times \boldsymbol{\omega}) \cdot \varphi
dxdt
$$
and
$$
\int_{\mathbb{R}^3}\boldsymbol{\omega}\cdot\varphi(\tau,x) dx -\int_{\mathbb{R}^3} \boldsymbol{\omega}_0\cdot\varphi(0,x) dx
$$
$$
=\int_0^{\tau}\int_{\mathbb{R}^3}\boldsymbol{\omega}\cdot\partial_t\varphi + \mathbf{v}\otimes\boldsymbol{\omega}:\nabla\varphi dxdt -\mu'\int_0^{\tau}\int_{\mathbb{R}^3}\nabla\boldsymbol{\omega}:\nabla\varphi dxdt
$$
$$
-(\mu'+\lambda')\int_0^{\tau}\int_{\mathbb{R}^3}\text{div}\boldsymbol{\omega}:\text{div}\varphi dxdt
+4\xi \int_0^{\tau}\int_{\mathbb{R}^3} \boldsymbol{\omega} \cdot \varphi
dxdt
-2\xi \int_0^{\tau}\int_{\mathbb{R}^3} \nabla \times \mathbf{v} dxdt,
$$
for any $\varphi\in C_c^{\infty}([0,T)\times\mathbb{R}^3)$, $\text{ div}\varphi = 0$.
We thus conclude the proof of the Theorem \ref{mainweak}.

\section{Weak--strong limit}

\subsection{Relative energy inequality}
\bigskip
Inspired by \cite{FeJiNo}, for any sufficiently smooth function $(r, \mathbf{U}, \mathbf{W})$, we introduce the following relative energy functional
$$
\mathcal{E}(\varrho_\varepsilon, \mathbf{u}_\varepsilon, \boldsymbol{\omega}_\varepsilon| r, \mathbf{U}, \mathbf{W})(\tau, x) 
$$   
\begin{equation} \label{entr-funct}
    = \int_{\mathbb{R}^3} \left( \frac{1}{2} \varrho_\varepsilon |\mathbf{u}_\varepsilon-\mathbf{U}|^2 + \frac{1}{2} \varrho_\varepsilon |\boldsymbol{\omega}_\varepsilon-\mathbf{W}|^2 + H(\varrho_\varepsilon)+H'(r)(\varrho_\varepsilon-r)-H(r) \right) dx. 
\end{equation}
Now, in order to derive a \textit{relative energy inequality} satisfied by the weak solution $(\varrho_\varepsilon, \mathbf{u}_\varepsilon, \boldsymbol{\omega}_\varepsilon)$, following \cite{Hu} we test the continuity equation (\ref{weak-cont}) with $\frac{1}{2}|\mathbf{U}|^2$, $\frac{1}{2}|\mathbf{W}|^2$ and $H'(r)-H'(1)$ respectively. We have,

$$
\int_{\mathbb{R}^3} \frac{1}{2} \varrho_\varepsilon(\tau, x) |\mathbf{U}|^2(\tau, x) dx - \int_{\mathbb{R}^3} \frac{1}{2 }\varrho_{0,\varepsilon} |\mathbf{U}|^2(0, x) dx 
$$   
\begin{equation} \label{test-U}   
    =  \int_0^T \int_{\mathbb{R}^3} \left( \varrho_\varepsilon \mathbf{U} \cdot \partial_t \mathbf{U} + \varrho_\varepsilon \mathbf{u}_\varepsilon \cdot \nabla \mathbf{U} \cdot \mathbf{U} \right) dxdt,
\end{equation}

$$
\int_{\mathbb{R}^3} \frac{1}{2} \varrho_\varepsilon(\tau, x) |\mathbf{W}|^2(\tau, x) dx - \int_{\mathbb{R}^3} \frac{1}{2 }\varrho_{0,\varepsilon} |\mathbf{W}|^2(0, x) dx 
$$   
\begin{equation} \label{test-W}
    =  \int_0^T \int_{\mathbb{R}^3} \left( \varrho_\varepsilon \mathbf{W} \cdot \partial_t \mathbf{W} + \varrho_\varepsilon \mathbf{u}_\varepsilon \cdot \nabla \mathbf{W} \cdot \mathbf{W} \right) dxdt
\end{equation}
and
$$
\int_{\mathbb{R}^3} \frac{1}{2} \varrho_\varepsilon(\tau, x) (H'(r)-H'(1))(\tau, x) dx - \int_{\mathbb{R}^3} \frac{1}{2 }\varrho_{0,\varepsilon} (H'(r)-H'(1))(0, x) dx 
$$   
\begin{equation} \label{test-H}
    =  \int_0^T \int_{\mathbb{R}^3} \left( \varrho_\varepsilon \partial_t H'(r) + \varrho_\varepsilon \mathbf{u}_\varepsilon \cdot \nabla H'(r) \right) dxdt.
\end{equation}
Now, we test (\ref{weak-mom}) and (\ref{weak-omega}) with $\mathbf{U}$ and $W$, respectively. We obtain

$$
\int_{\mathbb{R}^3} \varrho_\varepsilon \mathbf{u}_\varepsilon(\tau, x) \cdot \mathbf{U}(\tau, x) dx - \int_{\mathbb{R}^3} \varrho_{0,\varepsilon} \mathbf{u}_{0,\varepsilon} \cdot \mathbf{U}(0, x) dx
$$   
$$
= \int_0^T \int_{\mathbb{R}^3} \left( \varrho_\varepsilon \mathbf{u}_\varepsilon \cdot \partial_t \mathbf{U} + \varrho_\varepsilon \mathbf{u}_\varepsilon \otimes \mathbf{u}_\varepsilon : \nabla \mathbf{U} + \frac{p(\varrho_\varepsilon)}{\varepsilon^2} \text{div} \mathbf{U} \right) dxdt
$$
$$
- \nu \int_0^T \int_{\mathbb{R}^3} \left( (\mu + \xi) \nabla \mathbf{u}_\varepsilon : \nabla \mathbf{U} + (\mu + \lambda - \xi) \text{div} \mathbf{u}_\varepsilon \text{div} \mathbf{U} \right)dxdt
$$
\begin{equation} \label{test-mom}
+ \nu \int_0^T \int_{\mathbb{R}^3} 2\xi (\nabla \times \boldsymbol{\omega}_\varepsilon) \cdot \mathbf{U}) dxdt 
\end{equation}
and
$$
\int_{\mathbb{R}^3} \varrho_\varepsilon \boldsymbol{\omega}_\varepsilon (\tau, x) \cdot \mathbf{W}(\tau, x) dx - \int_{\mathbb{R}^3} \varrho_{0,\varepsilon} \boldsymbol{\omega}_{0,\varepsilon} \cdot \mathbf{W}(0, x) dx
$$   
$$
= \int_0^T \int_{\mathbb{R}^3} \left( \varrho_\varepsilon \boldsymbol{\omega}_\varepsilon \cdot \partial_t \mathbf{W} + \varrho_\varepsilon \boldsymbol{\omega}_\varepsilon \otimes \mathbf{u}_\varepsilon : \nabla \mathbf{W} - 4\xi \boldsymbol{\omega}_\varepsilon \cdot \mathbf{W} \right) dxdt
$$
$$
- \nu \int_0^T \int_{\mathbb{R}^3} \left( \mu' \nabla \boldsymbol{\omega}_\varepsilon : \nabla \mathbf{W} + (\mu' + \lambda') \text{div} \boldsymbol{\omega}_\varepsilon \text{div} \mathbf{W}
\right)dxdt
$$
\begin{equation} \label{test-omega}
+ \nu \int_0^T \int_{\mathbb{R}^3} 2\xi (\nabla \times \mathbf{u}_\varepsilon) \cdot \mathbf{W}) dxdt. 
\end{equation}
Plugging (\ref{test-U}) - (\ref{test-omega}) in (\ref{entr-funct}) and using the energy inequality (\ref{ei}), we obtain the following relative energy inequality
$$
\mathcal{E}(\varrho_\varepsilon, \mathbf{u}_\varepsilon, \boldsymbol{\omega}_\varepsilon| r, \mathbf{U}, \mathbf{W})(\tau, x)
- \mathcal{E}(\varrho_\varepsilon, \mathbf{u}_\varepsilon, \boldsymbol{\omega}_\varepsilon| r, \mathbf{U}, \mathbf{W})(0, x)
$$
$$
+\nu \int_0^T \int_{\mathbb{R}^3} \left( (\mu + \xi) \nabla\mathbf{u}_\varepsilon : (\nabla \mathbf{u}_\varepsilon -\nabla \mathbf{U}) + (\mu + \lambda -\xi) \text{div}\mathbf{u}_\varepsilon(\text{div}\mathbf{u}_\varepsilon-\text{div}\mathbf{U}) \right) dxdt
$$
$$
+\nu \int_0^T \int_{\mathbb{R}^3} \left( \mu' \nabla\boldsymbol{\omega}_\varepsilon : (\nabla \boldsymbol{\omega}_\varepsilon -\nabla \mathbf{W}) + (\mu' + \lambda') \text{div}\boldsymbol{\omega}_\varepsilon(\text{div}\boldsymbol{\omega}_\varepsilon-\text{div}\mathbf{W}) \right) dxdt
$$
$$
+\nu \int_0^T \int_{\mathbb{R}^3} \left( 4\xi \boldsymbol{\omega}_\varepsilon \cdot (\boldsymbol{\omega}_\varepsilon - \mathbf{W}) -2\xi (\nabla \times \boldsymbol{\omega}_\varepsilon) \cdot (\mathbf{u}_\varepsilon - \mathbf{U}) -2\xi (\nabla \times \mathbf{u}_\varepsilon) \cdot (\boldsymbol{\omega}_\varepsilon - \mathbf{W})
\right) dxdt
$$
$$
\leq
\int_0^T \int_{\mathbb{R}^3} \varrho_\varepsilon \left(  \partial_t \mathbf{U} +  \mathbf{u}_\varepsilon \cdot \nabla \mathbf{U}  \right) \cdot
\left( \mathbf{U} - \mathbf{u}_\varepsilon \right)dxdt
+
\int_0^T \int_{\mathbb{R}^3} \varrho_\varepsilon \left(  \partial_t \mathbf{W} +  \mathbf{u}_\varepsilon \cdot \nabla \mathbf{W}  \right) \cdot
\left( \mathbf{W} - \boldsymbol{\omega}_\varepsilon \right)dxdt
$$
$$
+\frac{1}{\varepsilon^2} \int_0^T \int_{\mathbb{R}^3} \left( (r-\varrho_\varepsilon)\partial_t H'(r) + \nabla H'(r) \cdot (r\mathbf{U}-\varrho_\varepsilon \mathbf{u}_\varepsilon) \right) dxdt
$$
\begin{equation} \label{rel-entr-ineq}
    - \frac{1}{\varepsilon^2} \int_0^T \int_{\mathbb{R}^3} \left( p(\varrho_\varepsilon) - p(r) \right) \text{div}\mathbf{U} dxdt,
\end{equation}
for any smooth function $(r, \mathbf{U}, \mathbf{W})$ such that
\begin{equation} \label{test_f}
r > 0, \ \ (r-1) \in C_c^\infty \left( [0,T] \times {\mathbb{R}^3} \right), \ \ (\mathbf{U}, \mathbf{W}) \in C_c^\infty \left( [0,T] \times {\mathbb{R}^3}; {\mathbb{R}^3} \right).
\end{equation}

\subsection{The incompressible inviscid limit}
After some algebra, the relative energy inequality (\ref{rel-entr-ineq}) could be rewritten as follows

$$
\mathcal{E}(\varrho_\varepsilon, \mathbf{u}_\varepsilon, \boldsymbol{\omega}_\varepsilon| r, \mathbf{U}, \mathbf{W})(\tau, \cdot)
- \mathcal{E}(\varrho_\varepsilon, \mathbf{u}_\varepsilon, \boldsymbol{\omega}_\varepsilon| r, \mathbf{U}, \mathbf{W})(0, \cdot)
$$
$$
+\nu \int_0^T \int_{\mathbb{R}^3} \left( (\mu + \xi) (\nabla\mathbf{u}_\varepsilon 
-\nabla \mathbf{U})
: (\nabla \mathbf{u}_\varepsilon -\nabla \mathbf{U}) + (\mu + \lambda -\xi) (\text{div}\mathbf{u}_\varepsilon -\text{div}\mathbf{U})(\text{div}\mathbf{u}_\varepsilon-\text{div}\mathbf{U}) \right) dxdt
$$
$$
+\nu \int_0^T \int_{\mathbb{R}^3} \left( \mu' (\nabla\boldsymbol{\omega}_\varepsilon -\nabla \mathbf{W}): (\nabla \boldsymbol{\omega}_\varepsilon -\nabla \mathbf{W}) + (\mu' + \lambda') (\text{div}\boldsymbol{\omega}_\varepsilon - \text{div}\mathbf{W})(\text{div}\boldsymbol{\omega}_\varepsilon-\text{div}\mathbf{W}) \right) dxdt
$$
$$
+\nu \int_0^T \int_{\mathbb{R}^3} \xi |2(\boldsymbol{\omega}_\varepsilon-\mathbf{W})-\nabla \times (\mathbf{u}_\varepsilon-\mathbf{U})|^2 dx
$$
$$
\leq
\int_0^T \int_{\mathbb{R}^3} \varrho_\varepsilon \left(  \partial_t \mathbf{U} +  \mathbf{u}_\varepsilon \cdot \nabla \mathbf{U}  \right) \cdot
\left( \mathbf{U} - \mathbf{u}_\varepsilon \right)dxdt
+
\int_0^T \int_{\mathbb{R}^3} \varrho_\varepsilon \left(  \partial_t \mathbf{W} +  \mathbf{u}_\varepsilon \cdot \nabla \mathbf{W}  \right) \cdot
\left( \mathbf{W} - \boldsymbol{\omega}_\varepsilon \right)dxdt
$$
$$
-\nu \int_0^T \int_{\mathbb{R}^3} \left( (\mu + \xi) \nabla \mathbf{U}
: (\nabla \mathbf{u}_\varepsilon -\nabla \mathbf{U}) + (\mu + \lambda -\xi) \text{div}\mathbf{U}(\text{div}\mathbf{u}_\varepsilon-\text{div}\mathbf{U}) \right) dxdt
$$
$$
-\nu \int_0^T \int_{\mathbb{R}^3} \left( \mu'  \nabla \mathbf{W}: (\nabla \boldsymbol{\omega}_\varepsilon -\nabla \mathbf{W}) + (\mu' + \lambda')   \text{div}\mathbf{W}(\text{div}\boldsymbol{\omega}_\varepsilon-\text{div}\mathbf{W}) \right) dxdt
$$
$$
-\nu \int_0^T \int_{\mathbb{R}^3} \left( 4\xi \mathbf{W} \cdot (\boldsymbol{\omega}_\varepsilon - \mathbf{W}) -2\xi (\nabla \times \mathbf{W}) \cdot (\mathbf{u}_\varepsilon - \mathbf{U}) -2\xi (\nabla \times \mathbf{U}) \cdot (\boldsymbol{\omega}_\varepsilon - \mathbf{W})
\right) dxdt
$$
$$
+\frac{1}{\varepsilon^2} \int_0^T \int_{\mathbb{R}^3} \left( (r-\varrho_\varepsilon)\partial_t H'(r) + \nabla H'(r) \cdot (r\mathbf{U}-\varrho_\varepsilon \mathbf{u}_\varepsilon) \right) dxdt
$$
\begin{equation} \label{rel-entr-ineq-2}
    - \frac{1}{\varepsilon^2} \int_0^T \int_{\mathbb{R}^3} \left( p(\varrho_\varepsilon) - p(r) \right) \text{div}\mathbf{U} dxdt.
\end{equation}
Now, let $r_{\varepsilon,\eta} = 1+ \varepsilon \psi_{\varepsilon,\eta}$. In the following, we will use
$$[{r}_{\varepsilon,\eta}, \mathbf{U}_{\varepsilon,\eta},\mathbf{W}_{\varepsilon,\eta}], \ \  \mathbf{U}_{\varepsilon,\eta}=(\mathbf{v} + \nabla \Psi_{\varepsilon,\eta}), \ \ \mathbf{W}_{\varepsilon,\eta} = \boldsymbol{\omega}$$
as test functions $\left[{r},\mathbf{U},\mathbf{W}\right]$ in the relative energy inequality (\ref{rel-entr-ineq-2}), with $(\mathbf{v}, \boldsymbol{\omega})$ strong solution of the system (\ref{cont-E}) - (\ref{micr-E}).

\subsubsection{Initial data}
For the initial data, similar to the analysis developed in \cite{Ca} and \cite{CaNe}, we have	
	$$
	\mathcal{E}_{\varepsilon,\eta} \left(\varrho_\varepsilon,\mathbf{u}_\varepsilon, \boldsymbol{\omega}_\varepsilon | {r}_{\varepsilon,\eta}, \mathbf{U}_{\varepsilon,\eta},\mathbf{W}_{\varepsilon,\eta} \right)\left(0, x \right)= \int_{\mathbb{R}^3}\frac{1}{2}\varrho_{0,\varepsilon}\left|\mathbf{u}_{0,\varepsilon}-\mathbf{u}_{0}\right|^{2} dx
	+\int_{\mathbb{R}^3}\frac{1}{2}\varrho_{0,\varepsilon}\left|\boldsymbol{\omega}_{0,\varepsilon}-\boldsymbol{\omega}_{0}\right|^{2} dx
	$$
	\begin{equation} \label{initial data conv}
	+\int_{\mathbb{R}^3}\frac{1}{\varepsilon^{2}}\left[H\left(1+\varepsilon\varrho_{0,\varepsilon}^{(1)}\right)-\varepsilon H^{\prime}\left(1+\varepsilon\varrho_{0}^{(1)}\right)\left(\varrho_{0,\varepsilon}^{(1)}-\varrho_{0}^{(1)}\right)-H\left(1+\varepsilon\varrho_{0}^{(1)}\right)\right]dx.
	\end{equation}
	For the first term on the right hand side of the equality (\ref{initial data conv}),
	we have
	$$
	\int_{\mathbb{R}^3}\frac{1}{2}\varrho_{0,\varepsilon}\left|\mathbf{u}_{0,\varepsilon}-\mathbf{u}_{0}\right|^{2}dx
	 =\int_{\mathbb{R}^3}\frac{1}{2}\left|1+\varepsilon\varrho_{0,\varepsilon}^{(1)}\right|\left|\mathbf{u}_{0,\varepsilon}-\mathbf{u}_{0}\right|^{2}dx
	$$
	$$
	 \leq\int_{\mathbb{R}^3}\frac{1}{2}\left|\mathbf{u}_{0,\varepsilon}-\mathbf{u}_{0}\right|^{2}dx+\int_{\mathbb{R}^3}\frac{1}{2}\left|\varepsilon\varrho_{0,\varepsilon}^{(1)}\right|\left|\mathbf{u}_{0,\varepsilon}-\mathbf{u}_{0}\right|^{2}dx
	$$
	$$
	\leq\int_{\mathbb{R}^3}\frac{1}{2}\left|\mathbf{u}_{0,\varepsilon}-\mathbf{u}_{0}\right|^{2}dx+\varepsilon\left\Vert \varrho_{0,\varepsilon}^{(1)}\right\Vert _{L^{\infty}\left(\mathbb{R}^{3}\right)}\int_{\mathbb{R}^3}\frac{1}{2}\left|\mathbf{u}_{0,\varepsilon}-\mathbf{u}_{0}\right|^{2}dx
	$$
	\begin{equation} \label{initial data conv1}
	\leq c\left(1+\varepsilon\right)\left\Vert \left|\mathbf{u}_{0,\varepsilon}-\mathbf{u}_{0}\right|^{2}\right\Vert _{L^{1}(\mathbb{R}^{3};\mathbb{R}^{3})}.
	\end{equation}
	Similar analysis can be done for the second term on the right hand side of (\ref{initial data conv}).
	For the third term on the right hand side of (\ref{initial data conv}),
	setting $a=1+\varepsilon\varrho_{0,\varepsilon}^{(1)}$ and $b=1+\varepsilon\varrho_{0}^{(1)}$
	and observing that
	$$
	H(a)=H(b)+H^{\prime}(b)(a-b)+\frac{1}{2}H^{\prime\prime}(\xi)(a-b)^{2},\;\;\xi\in\left(a,b\right),
	$$
	$$
	\left|H(a)-H^{\prime}(b)(a-b)-H(b)\right|\leq c\left|a-b\right|^{2},
	$$
	we have
	$$
	\int_{\mathbb{R}^3}\frac{1}{\varepsilon^{2}}\left[H\left(1+\varepsilon\varrho_{0,\varepsilon}^{(1)}\right)-\varepsilon H^{\prime}\left(1+\varepsilon\varrho_{0}^{(1)}\right)\left(\varrho_{0,\varepsilon}^{(1)}-\varrho_{0}^{(1)}\right)\right]dx
	$$
	$$
	\leq c\int_{\mathbb{R}^3}\frac{1}{\varepsilon^{2}}\left(\left|\varepsilon\left(\varrho_{0,\varepsilon}^{(1)}-\varrho_{0}^{(1)}\right)\right|^{2}\right)dx
	$$
	\begin{equation} \label{initial data conv2}
	\leq c\left\Vert \left|\varrho_{0,\varepsilon}^{(1)}-\varrho_{0}^{(1)}\right|^{2}\right\Vert _{L^{1}(\mathbb{R}^{3})}.
	\end{equation}
	We can conclude
	$$
	\mathcal{E}_{\varepsilon,\eta} \left(\varrho_\varepsilon,\mathbf{u}_\varepsilon, \boldsymbol{\omega}_\varepsilon | {r}_{\varepsilon,\eta}, \mathbf{U}_{\varepsilon,\eta},\mathbf{W}_{\varepsilon,\eta} \right)\left(0, x \right)
	$$
	\begin{equation} \label{conv-id}
	\leq c\left(1+\varepsilon\right)\left( \left\Vert \left|\mathbf{u}_{0,\varepsilon}-\mathbf{u}_{0}\right|^{2}\right\Vert _{L^{1}(\mathbb{R}^{3};\mathbb{R}^{3})}
	+\left\Vert \left|\boldsymbol{\omega}_{0,\varepsilon}-\boldsymbol{\omega}_{0}\right|^{2}\right\Vert _{L^{1}(\mathbb{R}^{3};\mathbb{R}^{3})}
	\right)
	+c\left\Vert \left|\varrho_{0,\varepsilon}^{(1)}-\varrho_{0}^{(1)}\right|^{2}\right\Vert _{L^{1}(\mathbb{R}^{3})}.
	\end{equation}

\subsubsection{The dissipative terms}\label{dpt}
We have
$$
\nu \int_0^T \int_{\mathbb{R}^3} (\mu+\xi) \nabla\mathbf{U}_{\varepsilon,\eta}: (\nabla\mathbf{U}_{\varepsilon,\eta}-\nabla\mathbf{u}_\varepsilon)dxdt
$$
$$
\leq \frac{\nu}{2}\int_0^\tau \int_{\mathbb{R}^3} (\mu+\xi) (\nabla\mathbf{U}_{\varepsilon,\eta}-\nabla\mathbf{u}_\varepsilon): (\nabla\mathbf{U}_{\varepsilon,\eta}-\nabla\mathbf{u}_\varepsilon)dxdt + \frac{\nu}{2} \int_0^\tau\int_{\mathbb{R}^3}(\mu+\xi) |\nabla\mathbf{U}_{\varepsilon,\eta}|^2dxdt,
$$
where we used the Young inequality. The first term can be absorbed by its counterpart on the left side of (\ref{rel-entr-ineq-2}) and the second term is estimated by $c(\eta,\mu, \xi, T)\nu$, which goes to zero as $\nu \to 0$. Similarly, is possible to treat the other viscous term and the micro-polar viscous terms. 

Moreover, thanks to the uniform bounds and the regularity of $\boldsymbol{\omega}$ and $\mathbf{U}_{\varepsilon,\eta}$, the term
$$
\nu \int_0^T \int_{\mathbb{R}^3} \left( 4\xi \boldsymbol{\omega} \cdot (\boldsymbol{\omega}_\varepsilon - \boldsymbol{\omega}) -2\xi (\nabla \times \boldsymbol{\omega}) \cdot (\mathbf{u}_\varepsilon - \mathbf{U}_{\varepsilon,\eta}) -2\xi (\nabla \times \mathbf{U}_{\varepsilon,\eta}) \cdot (\boldsymbol{\omega}_\varepsilon - \boldsymbol{\omega})
\right) dxdt
$$
goes to zero as $\nu \rightarrow 0$.

\subsubsection{The micro-polar convective term}\label{mct}
We have
$$
\int_0^T \int_{\mathbb{R}^3} \varrho_\varepsilon \left(  \partial_t \boldsymbol{\omega} +  \mathbf{u}_\varepsilon \cdot \nabla \boldsymbol{\omega} \right) \cdot
\left( \boldsymbol{\omega} - \boldsymbol{\omega}_\varepsilon \right)dxdt
$$
$$
= \int_0^T \int_{\mathbb{R}^3} \varrho_\varepsilon \partial_t \boldsymbol{\omega} \cdot (\boldsymbol{\omega} - \boldsymbol{\omega}_\varepsilon) dxdt
+ \int_0^T \int_{\mathbb{R}^3} \varrho_\varepsilon \mathbf{u}_\varepsilon \cdot \nabla \boldsymbol{\omega} \cdot (\boldsymbol{\omega} - \boldsymbol{\omega}_\varepsilon) dxdt
$$
Using (\ref{micr-E}), we have
$$
\int_0^T \int_{\mathbb{R}^3} \partial_t \boldsymbol{\omega} \cdot \varrho_\varepsilon (\boldsymbol{\omega} - \boldsymbol{\omega}_\varepsilon) dxdt
= - \int_0^T \int_{\mathbb{R}^3} \mathbf{v} \cdot \nabla \boldsymbol{\omega} \cdot \varrho_\varepsilon (\boldsymbol{\omega} - \boldsymbol{\omega}_\varepsilon) dxdt
$$
$$
= - \int_0^T \int_{\mathbb{R}^3} \varrho_\varepsilon \mathbf{U}_{\varepsilon,\eta} \cdot \nabla{\boldsymbol{\omega}} \cdot (\boldsymbol{\omega} - \boldsymbol{\omega}_\varepsilon) dxdt  + \int_0^T \int_{\mathbb{R}^3} \varrho_\varepsilon \nabla \Psi_{\varepsilon,\eta} \cdot \nabla \boldsymbol{\omega} \cdot (\boldsymbol{\omega} - \boldsymbol{\omega}_\varepsilon) dxdt.
$$
Consequently, after some algebra, we have
$$
\int_0^T \int_{\mathbb{R}^3} \varrho_\varepsilon \left(  \partial_t \boldsymbol{\omega} +  \mathbf{u}_\varepsilon \cdot \nabla \boldsymbol{\omega}  \right) \cdot
\left( \boldsymbol{\omega} - \boldsymbol{\omega}_\varepsilon \right)dxdt
$$
$$
= \int_0^T \int_{\mathbb{R}^3} \varrho_\varepsilon (\mathbf{u}_\varepsilon - \mathbf{U}_{\varepsilon,\eta}) \cdot \nabla \boldsymbol{\omega} \cdot (\boldsymbol{\omega} - \boldsymbol{\omega}_\varepsilon) dxdt
+ \int_0^T \int_{\mathbb{R}^3} \varrho_\varepsilon \nabla \Psi_{\varepsilon,\eta} \cdot \nabla \boldsymbol{\omega} \cdot (\boldsymbol{\omega} - \boldsymbol{\omega}_\varepsilon) dxdt.
$$
We can estimate the first term as follows
$$
\int_0^T \int_{\mathbb{R}^3} \varrho_\varepsilon (\mathbf{u}_\varepsilon - \mathbf{U}_{\varepsilon,\eta}) \cdot \nabla \boldsymbol{\omega} \cdot (\boldsymbol{\omega} - \boldsymbol{\omega}_\varepsilon) dxdt
$$
$$
\leq c(\eta,T) 
\left(
\int_0^T \int_{\mathbb{R}^3} \varrho_\varepsilon |\mathbf{u}_\varepsilon - \mathbf{U}_{\varepsilon,\eta}|^2  dxdt
+\int_0^T \int_{\mathbb{R}^3} \varrho_\varepsilon |\boldsymbol{\omega} - \boldsymbol{\omega}_\varepsilon|^2  dxdt
\right)
$$
\begin{equation} \label{est-micr}
    \leq c(\eta,T) \int_0^T \mathcal{E}(t,\cdot) dx,
\end{equation}
where we used the Young inequality and the regularity of $\boldsymbol{\omega}$. 
While, for the acoustic contribution we have
$$
\int_0^T \int_{\mathbb{R}^3} \varrho_\varepsilon \nabla \Psi_{\varepsilon,\eta} \cdot \nabla \boldsymbol{\omega} \cdot (\boldsymbol{\omega} - \boldsymbol{\omega}_\varepsilon) dxdt
$$
$$
= \int_0^T \int_{\mathbb{R}^3} \varrho_\varepsilon \nabla \Psi_{\varepsilon,\eta} \cdot \nabla \boldsymbol{\omega} \cdot \boldsymbol{\omega} dxdt
- \int_0^T \int_{\mathbb{R}^3} \varrho_\varepsilon \boldsymbol{\omega}_\varepsilon \cdot \nabla \Psi_{\varepsilon,\eta} \cdot \nabla \boldsymbol{\omega} dxdt.
$$
$$
= \varepsilon \int_0^T \int_{\mathbb{R}^3} \frac{\varrho_\varepsilon - 1}{\varepsilon} \nabla \Psi_{\varepsilon,\eta} \cdot \nabla \boldsymbol{\omega} \cdot \boldsymbol{\omega} dxdt
+\int_0^T \int_{\mathbb{R}^3} \nabla \Psi_{\varepsilon,\eta} \cdot \nabla \boldsymbol{\omega} \cdot \boldsymbol{\omega} dxdt
$$
$$
- \int_0^T \int_{\mathbb{R}^3} \varrho_\varepsilon \boldsymbol{\omega}_\varepsilon \cdot \nabla \Psi_{\varepsilon,\eta} \cdot \nabla \boldsymbol{\omega} dxdt.
$$
According to (\ref{unif_bound1}) - (\ref{unif_bound2}), the first term can be estimated as follows
$$
c(T)\left\|\frac{\varrho_\varepsilon-1}{\varepsilon}\right\|_{L^\infty_T(L^2(\mathbb{R}^3)+L^{\gamma}(\mathbb{R}^3))} 
\left\|\nabla \boldsymbol{\omega} \right\|_{L^\infty_T(L^8(\mathbb{R}^3)+L^{\frac{8\gamma}{3\gamma-4}}(\mathbb{R}^3))} 
$$
$$
\cdot \left\| \boldsymbol{\omega} \right\|_{L^\infty_T(L^8(\mathbb{R}^3)+L^{\frac{8\gamma}{3\gamma-4}}(\mathbb{R}^3))} 
\left\|\nabla\Psi_{\varepsilon,\eta}\right\|_{L^{4}_T(L^4(\mathbb{R}^3)+L^{4}(\mathbb{R}^3))}
$$
\begin{equation}\label{est-micr-1}
\leq c(\eta,T)\varepsilon^{\frac{1}{4}},
\end{equation}
as well as the second term,
\begin{equation}\label{est-micr-2}
    c(\eta,T)\varepsilon^{\frac{1}{4}},
\end{equation}
thanks to the dispersive estimates and the regularity of $\boldsymbol{\omega}$.
Finally, for the last term we have the following estimate
$$
c(\eta,T) \left\| \varrho_\varepsilon \boldsymbol{\omega}_\varepsilon\right\|_{L^\infty_T(L^2(\mathbb{R}^3)+L^{\frac{2\gamma}{\gamma+1}}(\mathbb{R}^3))} \left\|\nabla \boldsymbol{\omega} \right\|_{L^\infty_T(L^{4}(\mathbb{R}^3)+L^\infty(\mathbb{R}^3))}
$$
$$
\cdot \left\|\nabla\Psi_{\varepsilon,\eta}\right\|_{L^{4}_T(L^{4}(\mathbb{R}^3))+L^{2\gamma}_T(L^{\frac{2\gamma}{\gamma -1}}(\mathbb{R}^3))}
$$
\begin{equation}\label{est-micr-3}
\leq c(\eta,T) \varepsilon^{\min\{ \frac{1}{4}, \frac{1}{2\gamma} \}}.
\end{equation}
In conclusion, the convective term related to the momentum equation and the pressure terms are treated as in \cite{Ca} and \cite{CaNe} (see the Appendix).   

\subsection{Proof of Theorem \ref{maineuler}}
From the above estimates, (\ref{r1final}), (\ref{prfinal}) and the conservation of the acoustic energy (\ref{acwenergy}), we obtain
$$
\mathcal{E}_{\varepsilon,\eta}\left(\varrho,\mathbf{u}, \boldsymbol{\omega} \mid {r},\mathbf{U},W\right)\left(\tau, x \right)
\leq c(\eta,T)\varepsilon^{\alpha} + \int_0^T c(t) \mathcal{E}_{\varepsilon,\eta}(t) dt, 
\ \ \alpha=\min\{\frac{1}{4},\frac{1}{2\gamma}\}.
$$
By Gronwall's inequality, we have
\begin{equation}\label{cpt1}
\mathcal{E}_{\varepsilon,\eta}\left(\varrho,\mathbf{u}, \boldsymbol{\omega} \mid {r},\mathbf{U},W\right)\left(\tau, x \right)
\leq c(\eta,T)\varepsilon^\alpha + c(T) \mathcal{E}_{\varepsilon,\eta}(0, x), \ \ \text{ a.e. } \tau\in (0,T).
\end{equation}
The convergence rate (\ref{conv-r}) given in Corollary \ref{cor} comes from (\ref{conv-id}) together with (\ref{cpt1}). 
Now, sending $\varepsilon\to 0$,
according to the assumptions (\ref{comp}), 
we find
$$
\lim_{\varepsilon\to 0}\mathcal{E}\left(\varrho_{\varepsilon},\mathbf{u}_{\varepsilon}, \boldsymbol{\omega}_{\varepsilon}\mid {{r}_{\varepsilon,\eta}},\mathbf{U_{\varepsilon,\eta}}, W_{ \varepsilon,\eta}\right)\left(\tau, x \right)=0 \text{ uniformly in }\tau\in (0,T),
$$
where $r_{\varepsilon,\eta} = 1+ \psi_{\varepsilon,\eta},\,\mathbf{U}_{\varepsilon,\eta}=\mathbf{v} + \nabla \Psi_{\varepsilon,\eta}, \, W_{\varepsilon,\eta} = \boldsymbol{\omega}$. We thus conclude the proof of Theorem \ref{maineuler}.
Indeed, $\nabla \Psi_{\varepsilon,\eta}\to 0$ in $L^q(0,T;L^p(\mathbb{R}^3))$ as $\varepsilon\to 0$ for any $(p,q)>2$ according to (\ref{acw_6}). Consequently, for any compact set $K\subset\mathbb{R}^3$, we have
$$
\left\|\sqrt{\varrho_{\varepsilon}}\mathbf{u}_{\varepsilon} - \mathbf{v}\right\|_{L^2_T(L^2(K))} \leq \left\|\sqrt{\varrho_{\varepsilon}}\mathbf{u}_{\varepsilon}
 - \mathbf{U}_{\varepsilon,\eta}\right\|_{L^2_T(L^2(\mathbb{R}^3))}
$$
$$
+ c(T,K)\left\|{\nabla{\Psi}_{\varepsilon,\eta}} \right\|_{L^q_T(L^p(K))},
$$
which goes to zero as $\varepsilon\to 0$.

\appendix
\section{Appendix}

For the reader's convenience, we recall the analysis from \cite{Ca} and \cite{CaNe}
providing the estimates on the convective and pressure terms.

\subsection{The convective terms} We decompose
$$
\int_0^T \int_{\mathbb{R}^3} \varrho_\varepsilon \left(  \partial_t \mathbf{U}_{\varepsilon,\eta} +  \mathbf{u_\varepsilon} \cdot \nabla \mathbf{U}_{\varepsilon,\eta}  \right) \cdot
\left( \mathbf{U}_{\varepsilon,\eta} - \mathbf{u}_\varepsilon \right)dxdt
$$
$$
= 
\int_0^T\int_{\mathbb{R}^3}\varrho_\varepsilon\left(\partial_{t}\mathbf{U}_{\varepsilon,\eta}+\mathbf{U}_{\varepsilon,\eta}\cdot\nabla\mathbf{U}_{\varepsilon,\eta}\right)\cdot\left(\mathbf{U}_{\varepsilon,\eta}-\mathbf{u}_\varepsilon\right)dxdt
$$
\begin{equation}\label{re1_1}
-\int_0^T \int_{\mathbb{R}^3} \varrho_\varepsilon\left(\mathbf{U}_{\varepsilon,\eta}-\mathbf{u}_\varepsilon\right)\cdot\nabla\mathbf{U}_{\varepsilon,\eta}\cdot\left(\mathbf{U}_{\varepsilon,\eta}-\mathbf{u}_\varepsilon\right) dxdt.
\end{equation}
The last term is estimated as follows
\begin{equation} \label{re1_2}
\int_0^T \int_{\mathbb{R}^3} \varrho_\varepsilon\left(\mathbf{U}_{\varepsilon,\eta}-\mathbf{u}_\varepsilon\right)\cdot\nabla\mathbf{U}_{\varepsilon,\eta}\cdot\left(\mathbf{U}_{\varepsilon,\eta}-\mathbf{u}_\varepsilon\right) dxdt
\leq
\int_0^T c(t)\mathcal{E}(t, \cdot) dt.
\end{equation}
For the first term on the right side of (\ref{re1_1}), we have
$$
\int_0^T \int_{\mathbb{R}^3}\varrho_\varepsilon\left(\partial_{t}\mathbf{U}_{\varepsilon,\eta}+\mathbf{U}_{\varepsilon,\eta}\cdot\nabla\mathbf{U}_{\varepsilon,\eta}\right)\cdot\left(\mathbf{U}_{\varepsilon,\eta}-\mathbf{u}_\varepsilon\right) dxdt
$$
$$
=\int_0^T \int_{\mathbb{R}^3}\varrho_\varepsilon\left(\partial_{t}\mathbf{v}+\mathbf{v}\cdot\nabla \mathbf{v}\right)\cdot\left(\mathbf{U}_{\varepsilon,\eta}-\mathbf{u}_\varepsilon\right)dxdt
+ \int_0^T \int_{\mathbb{R}^3}\varrho_\varepsilon \partial_{t}\nabla {\Psi_{\varepsilon,\eta}}\cdot\left(\mathbf{U}_{\varepsilon,\eta}-\mathbf{u_\varepsilon}\right) dxdt
$$
$$
+\int_0^T \int_{\mathbb{R}^3}\varrho_\varepsilon\nabla {\Psi_{\varepsilon,\eta}}\cdot\nabla \nabla{\Psi_{\varepsilon,\eta}}\cdot\left(\mathbf{U}_{\varepsilon,\eta}-\mathbf{u}_\varepsilon\right)dxdt
$$
\begin{equation} \label{re3}
+\int_0^T \int_{\mathbb{R}^3} \varrho_\varepsilon\left(\mathbf{v}\cdot\nabla (\nabla \Psi_{\varepsilon,\eta}) + \nabla \Psi_{\varepsilon,\eta}\cdot\nabla \mathbf{v} \right)\cdot\left(\mathbf{U}_{\varepsilon,\eta}-\mathbf{u}_\varepsilon\right)dxdt.
\end{equation}
Next,
$$
\int_0^T \int_{\mathbb{R}^3}\varrho_\varepsilon\left(\partial_{t}\mathbf{v}+\mathbf{v}\cdot\nabla\mathbf{v}\right)\cdot\left(\mathbf{U}_{\varepsilon,\eta}-\mathbf{u}_\varepsilon\right)dxdt=
{I_{1}+I_{2}}.
$$
We have
$$
{I_{1}}=\int_0^T \int_{\mathbb{R}^3}  \varrho_\varepsilon \mathbf{u}_\varepsilon \cdot \nabla \Pi dxdt 
$$
$$
= \int_{\mathbb{R}^3}  \varrho_\varepsilon  \Pi dx\left|_{t=0}^T \right.
- \int_0^T \int_{\mathbb{R}^3}  \varrho_\varepsilon \partial_t \Pi dxdt
$$
\begin{equation} \label{d3}
=\varepsilon \int_{\mathbb{R}^3} \frac{\varrho_\varepsilon - 1}{\varepsilon} \Pi dx\left|_{t=0}^T \right.
- \varepsilon \int_0^T \int_{\mathbb{R}^3}  \frac{\varrho_\varepsilon -1}{\varepsilon} \partial_t \Pi dxdt \leq c(\eta,T)\varepsilon
\end{equation}
according to (\ref{unif_bound1}) - (\ref{unif_bound2}), and
	$$
	|{I_{2}}|=\left|\int_{0}^{T}\int_{\mathbb{R}^3}
	\varrho_\varepsilon\mathbf{U}_{\varepsilon,\eta}\cdot\nabla \Pi dxdt\right|\leq\left|\int_{0}^{T}\int_{\mathbb{R}^3}
	\left(\varrho_\varepsilon-1\right)\cdot\mathbf{U}_{\varepsilon,\eta}\cdot\nabla \Pi dxdt\right|
	$$
	\begin{equation} \label{split}
	+\left|\int_{0}^{T}\int_{\mathbb{R}^3}
	\mathbf{U}_{\varepsilon,\eta}\cdot\nabla \Pi dxdt \right|.
	\end{equation}
For the first term on the right hand side of (\ref{split}), we have
$$
\left|\int_{0}^{T}\int_{\mathbb{R}^3}\left(\varrho_\varepsilon-1\right)\cdot\mathbf{U}_{\varepsilon,\eta}\cdot\nabla\Pi dxdt \right|
\leq \varepsilon\left|\int_{0}^{T}\int_{\mathbb{R}^3}\frac{\left(\varrho_\varepsilon-1\right)}{\varepsilon}\cdot\mathbf{U}_{\varepsilon,\eta}\cdot\nabla\Pi dxdt \right|
$$
$$
\leq c(T)\varepsilon,
$$
according to (\ref{unif_bound1}) - (\ref{unif_bound2}) and the energy estimate (\ref{acw_3}). For the second term on the right hand side of (\ref{split}), we have
\begin{equation} \label{press_conv2-4}
\int_{0}^{T} \int_{\mathbb{R}^3}\mathbf{U}_{\varepsilon,\eta}\cdot\nabla \Pi dxdt =
\int_{0}^{T}\int_{\mathbb{R}^3}\mathbf{v}\cdot\nabla \Pi dxdt+\int_{0}^{T} \int_{\mathbb{R}^3}\nabla \Psi_{\varepsilon,\eta}\cdot\nabla \Pi dxdt. 
\end{equation}
Performing integration by parts in the first term on the right-hand
side of (\ref{press_conv2-4}), we have
$$
\int_{0}^{T}\int_{\mathbb{R}^3}\textrm{div}\mathbf{v}\cdot\Pi dxdt=0,
$$
thanks to the incompressibility condition $\textrm{div}\mathbf{v}=0$.
For the second term on the right-hand side of (\ref{press_conv2-4}),
using integration by parts and acoustic equation, we
have
$$
\int_{0}^{T}\int_{\mathbb{R}^3}\nabla\Psi_{\varepsilon,\eta}\cdot\nabla\Pi dxdt =-\int_{0}^{T}\int_{\mathbb{R}^3}\Delta\Psi_{\varepsilon,\eta}\cdot\Pi dxdt
$$
$$
=\varepsilon\int_{0}^{T}\int_{\mathbb{R}^3}\partial_{t}\psi_{\varepsilon,\eta}\cdot\Pi dxdt
$$
\begin{equation} \label{phi_p}
=\varepsilon\left[\int_{\mathbb{R}^3}\psi_{\varepsilon,\eta}\cdot\Pi dx \right]_{t=0}^{t=T}-\varepsilon\int_{0}^{T}\int_{\mathbb{R}^3}\psi_{\varepsilon,\eta}\cdot\partial_{t}\Pi dxdt,
\end{equation}
that goes to zero for $\varepsilon\rightarrow0$.\\

Moreover, by using similar argument as above,
the last two terms in (\ref{re3}) are of order
\begin{equation}\label{re5}
c(\eta,T)(1+\varepsilon) \|\nabla\Psi_{\varepsilon,\eta} \|_{L^4_T(W^{1,4}(\mathbb{R}^3;\mathbb{R}^3))} \leq c(\eta,T)\varepsilon^{\frac{1}{4}}.
\end{equation}
Now, using ${\rm div}\mathbf{v}=0$, we have
$$
\int_0^T\int_{\mathbb{R}^3}\varrho_\varepsilon \partial_{t}\nabla{\Psi_{\varepsilon,\eta}}\cdot\left(\mathbf{U}_{\varepsilon,\eta}-\mathbf{u}_\varepsilon\right) dxdt = - \int_0^T\int_{\mathbb{R}^3}\varrho_\varepsilon\mathbf{u}_\varepsilon\cdot \partial_{t}\nabla{\Psi_{\varepsilon,\eta}} dxdt
$$
\begin{equation} \label{re6}
+ \int_0^T \int_{\mathbb{R}^3}(\varrho_\varepsilon - 1) \mathbf{v}\cdot\partial_{t}\nabla\Psi_{\varepsilon,\eta} dxdt + \int_0^T\int_{\mathbb{R}^3}\varrho_\varepsilon \partial_{t} \nabla\Psi_{\varepsilon,\eta} \cdot \nabla\Psi_{\varepsilon,\eta} dxdt
\end{equation}
The first term on the right side of (\ref{re6}) will be cancelled later by the pressure term while, by using the acoustic wave equations (\ref{acw_1}),
the second term equals to
$$
\int_0^T
\int_{\mathbb{R}^3}
\frac{\varrho_\varepsilon-1}{\varepsilon} \varepsilon \partial_{t} \nabla {\Psi_{\varepsilon,\eta}}\cdot \mathbf{v} dxdt = 
-\int_0^T \int_{\mathbb{R}^3}
\frac{\varrho_\varepsilon-1}{\varepsilon} a^2\nabla\psi_{\varepsilon,\eta}\cdot \mathbf{v} dxdt
$$
$$
\leq c(T)\left\|\frac{\varrho_\varepsilon-1}{\varepsilon}\right\|_{L^\infty_T(L^2(\mathbb{R}^3)+L^{\gamma}(\mathbb{R}^3))} \left\|\mathbf{v}\right\|_{L^\infty_T(L^4(\mathbb{R}^3;\mathbb{R}^3)+L^{\frac{4\gamma}{3\gamma-4}}(\mathbb{R}^3;\mathbb{R}^3))} \left\|\nabla\psi_{\varepsilon,\eta}\right\|_{L^{4}_T(L^4(\mathbb{R}^3)+L^{4}(\mathbb{R}^3))}
$$
\begin{equation}\label{re7}
 \leq c(\eta,T)\varepsilon^{\frac{1}{4}},
\end{equation}
according to (\ref{unif_bound1}) - (\ref{unif_bound2}). Finally, by using again the acoustic equations $\varepsilon\partial_t\nabla\Psi_{\varepsilon,\eta} = -a^{2} \nabla\psi_{\varepsilon,\eta}$,
$$
\int_0^T \int_{\mathbb{R}^3}\varrho_\varepsilon \partial_{t} \nabla\Psi_{\varepsilon,\eta} \cdot \nabla\Psi_{\varepsilon,\eta} dxdt
$$
$$
= - a^{2}\int_0^T \int_{\mathbb{R}^3}\frac{\varrho_\varepsilon-1}{\varepsilon}  \nabla \psi_{\varepsilon,\eta} \cdot \nabla \Psi_{\varepsilon,\eta} dxdt + \frac{1}{2}\int_{\mathbb{R}^3}|\nabla \Psi_{\varepsilon,\eta}|^2 dx\left|_{t=0}^{T}\right.
$$
\begin{equation}\label{re9}
\leq c(\eta,T)\varepsilon^{\frac{1}{4}} + \frac{1}{2}\int_{\mathbb{R}^3}|\nabla \Psi_{\varepsilon,\eta}|^2 dx\left|_{t=0}^{T}\right..
\end{equation}
From (\ref{re1_1}) to (\ref{re9}) we conclude that
$$
\int_0^T \int_{\mathbb{R}^3} \varrho_\varepsilon \left(  \partial_t \mathbf{U}_{\varepsilon,\eta} +  \mathbf{u}_\varepsilon \cdot \nabla \mathbf{U}_{\varepsilon,\eta}  \right) \cdot
\left( \mathbf{U}_{\varepsilon,\eta} - \mathbf{u}_\varepsilon \right)dxdt \leq c(\eta,T)\varepsilon^{\frac{1}{4}} + \int_0^T c(t)\mathcal{E}(t, \cdot) dt
$$
\begin{equation}\label{r1final}
+ \frac{1}{2}\int_{\mathbb{R}^3}|\nabla\Psi_{\varepsilon,\eta}|^2 dx\left|_{t=0}^{T
}\right. - \int_0^T \int_{\mathbb{R}^3}\varrho_\varepsilon\mathbf{u}_\varepsilon\cdot \partial_{t}\nabla{\Psi_{\varepsilon,\eta}} dxdt.
\end{equation}

\subsection{The pressure terms}
We recall
$$
\frac{1}{\varepsilon^{2}}\int_0^T \int_{\mathbb{R}^3}\left(r-\varrho_\varepsilon\right)\partial_{t}H^{\prime}\left(r\right)
+ \nabla H^{\prime}\left(r\right)\cdot(r \mathbf{U}_{\varepsilon,\eta}-\varrho_\varepsilon\mathbf{u}_\varepsilon) dxdt
$$
$$
- \frac{1}{\varepsilon^{2}}\int_0^T \int_{\mathbb{R}^3}
(p\left(\varrho_\varepsilon \right)-p\left(r \right))\textrm{div}\mathbf{U}_{\varepsilon,\eta} 
dxdt
$$
where $r=r_{\varepsilon,\eta} = 1 + \varepsilon \psi_{\varepsilon,\eta}$. We have,
$$
\int_0^T \int_{\mathbb{R}^3} 
\nabla H^{\prime}\left(r\right) \cdot r \mathbf{U}_{\varepsilon,\eta}
dxdt
= - \int_0^T \int_{\mathbb{R}^3} p(r)\text{div}\mathbf{U}_{\varepsilon,\eta} dxdt
$$
that cancels with its counterpart.
Next,
$$
\frac{1}{\varepsilon^{2}}\int_0^T \int_{\mathbb{R}^3} \varrho_\varepsilon\mathbf{u}_\varepsilon\cdot \nabla H^{\prime}\left(r\right) dx =\frac{1}{\varepsilon}\int_0^T \int_{\mathbb{R}^3}  \varrho_\varepsilon\mathbf{u}_\varepsilon\cdot \nabla\psi_{\varepsilon,\eta}  {H}^{\prime\prime}(r) dxdt
$$
$$
=\int_0^T \int_{\mathbb{R}^3} \varrho_\varepsilon\mathbf{u}_\varepsilon\cdot \nabla\psi_{\varepsilon,\eta} \frac{ {H}^{\prime\prime}(1+\varepsilon\psi_{\varepsilon,\eta}) - {H}^{\prime\prime}(1)}{\varepsilon} dxdt + \frac{1}{\varepsilon}\int_0^T \int_{\mathbb{R}^3} a^2 \varrho_\varepsilon\mathbf{u}_\varepsilon\cdot \nabla\psi_{\varepsilon,\eta}  dxdt
$$
where ${H}^{\prime\prime}(1)=p'(1)=a^2$. Observing  that
$$\left| \frac{ {H}^{\prime\prime}(1+\varepsilon\psi_{\varepsilon,\eta}) - {H}^{\prime\prime}(1)}{\varepsilon}\right|\leq c | \psi_{\varepsilon,\eta} |, $$
the first term on the right side can be estimated in the following way
$$
c(\eta,T) \left\| \varrho_\varepsilon\mathbf{u}_\varepsilon\right\|_{L^\infty_T(L^2(\mathbb{R}^3;\mathbb{R}^3)+L^{\frac{2\gamma}{\gamma+1}}(\mathbb{R}^3;\mathbb{R}^3))} \left\|\psi\right\|_{L^\infty_T(L^{4}(\mathbb{R}^3)+L^\infty(\mathbb{R}^3))} \left\|\nabla\psi\right\|_{L^{4}_T(L^{4}(\mathbb{R}^3))+L^{2\gamma}_T(L^{\frac{2\gamma}{\gamma -1}}(\mathbb{R}^3))}
$$
\begin{equation}\label{pr1}
\leq c(\eta,T) \varepsilon^{\min\{ \frac{1}{4}, \frac{1}{2\gamma} \}}
\end{equation}
For the second term, by using the acoustic equations,
$$
\frac{1}{\varepsilon}\int_0^T \int_{\mathbb{R}^3} a^2 \varrho_\varepsilon\mathbf{u}_\varepsilon\cdot \nabla\psi_{\varepsilon,\eta}  dxdt
= - \int_0^\tau\int_{\mathbb{R}^3}  \varrho_\varepsilon\mathbf{u}_\varepsilon\cdot \partial_t\Psi_{\varepsilon,\eta}  dxdt,
$$
which cancels with its counterpart on the right side of (\ref{r1final}) (or (\ref{re6})). Now, we have
$$
\frac{1}{\varepsilon^{2}}\int_0^T \int_{\mathbb{R}^3}\left(r-\varrho_\varepsilon\right)\partial_{t}H^{\prime}\left(r\right) - p\left(\varrho_\varepsilon \right)\textrm{div}\mathbf{U}_{\varepsilon,\eta}  dxdt
$$
$$
= \int_0^T \int_{\mathbb{R}^3} \frac{1-\varrho_\varepsilon}{\varepsilon} H^{\prime\prime}(r)\partial_t\psi_{\varepsilon,\eta} dxdt + \int_0^T \int_{\mathbb{R}^3} \psi_{\varepsilon,\eta} H^{\prime\prime}(r)\partial_t\psi_{\varepsilon,\eta} dx_hdt
$$
$$
- \int_0^T \int_{\mathbb{R}^3}\frac{ p (\varrho_\varepsilon)-p'(1)(\varrho_\varepsilon-1)-p(1)}{\varepsilon^2} \Delta {\Psi_{\varepsilon,\eta}} dxdt 
$$
\begin{equation}\label{pr3}
- \int_0^T \int_{\mathbb{R}^3}  p'(1)\frac{\varrho_\varepsilon -1}{\varepsilon^2}\Delta\Psi_{\varepsilon,\eta} dxdt.
\end{equation}
Then, the following estimate holds
\begin{equation}\label{pr8}
\int_0^T \int_{\mathbb{R}^3}\frac{ p (\varrho_\varepsilon)-p'(1)(\varrho_\varepsilon-1)-p(1)}{\varepsilon^2} \Delta {\Psi_{\varepsilon,\eta}} dxdt \leq c(\eta,T)\varepsilon^{\frac{1}{4}}.
\end{equation}
Now, we write
$$
\int_0^T \int_{\mathbb{R}^3} \frac{ 1-\varrho_\varepsilon}{\varepsilon} H^{\prime\prime}(r)\partial_t\psi_{\varepsilon,\eta} dxdt = \int_0^T \int_{\mathbb{R}^3} \frac{1- \varrho_\varepsilon}{\varepsilon} H^{\prime\prime}(1)\partial_t\psi_{\varepsilon,\eta} dxdt
$$
$$
+ \int_0^T \int_{\mathbb{R}^3} \frac{ 1-\varrho_\varepsilon}{\varepsilon} \left(H^{\prime\prime}(r)-H^{\prime\prime}(1)\right)\partial_t\psi_{\varepsilon,\eta} dxdt .
$$
Using the acoustic equation, the first term on the right side cancels with the last term in (\ref{pr3}) while the remaining term equals to
$$
- \int_0^T \int_{\mathbb{R}^3} \frac{1 - \varrho_\varepsilon}{\varepsilon} \frac{H^{\prime\prime}(r)-H^{\prime\prime}(1)}{\varepsilon}\Delta\Psi_{\varepsilon,\eta} dxdt
$$
and can be estimated as follows
$$
c(T)\left\|\frac{\varrho_\varepsilon -1}{\varepsilon}\right\|_{L^\infty_T(L^{2}(\mathbb{R}^3)+L^{\gamma}(\mathbb{R}^3))} \left\|\psi_{\varepsilon,\eta}\right\|_{L^\infty_T(L^{4}(\mathbb{R}^3)+L^{\frac{4\gamma}{3\gamma -4}}(\mathbb{R}^3))} \left\|\Delta\Psi_{\varepsilon,\eta}\right\|_{L^{4}_T(L^4(\mathbb{R}^3;\mathbb{R}^3)+L^4(\mathbb{R}^3;\mathbb{R}^3))}
$$
\begin{equation}\label{pr4}
\leq c(\eta,T)\varepsilon^{\frac{1}{4}}.
\end{equation}
Similarly,
$$
\int_0^T \int_{\mathbb{R}^3}  \psi_{\varepsilon,\eta} H^{\prime\prime}(r)\partial_t\psi_{\varepsilon,\eta} dxdt =  \int_0^T \int_{\mathbb{R}^3} \psi_{\varepsilon,\eta} H^{\prime\prime}(1) \partial_t\psi_{\varepsilon,\eta} dxdt
$$
$$
+ \int_0^T \int_{\mathbb{R}^3}  \psi_{\varepsilon,\eta} \left( H^{\prime\prime}(r)-H^{\prime\prime}(1)\right)\partial_t\psi_{\varepsilon,\eta} dxdt
$$
$$
\leq \frac{1}{2} \int_{\mathbb{R}^3} a^2 \left|\psi_{\varepsilon,\eta}\right|^2 dx \left|_{t=0}^T \right. + c(T)\left\|\psi_{\varepsilon,\eta}\right\|_{L^\infty_T(L^2(\mathbb{R}^3))} \left\|\psi_{\varepsilon,\eta}\right\|_{L^\infty_T(L^4(\mathbb{R}^3))} \left\|\Delta\Psi_{\varepsilon,\eta}\right\|_{L^{4}_T(L^4(\mathbb{R}^3;\mathbb{R}^3))}
$$
\begin{equation}\label{pr6}
\leq \frac{1}{2} \int_{\mathbb{R}^3} a^2 \left|\psi_{\varepsilon,\eta}\right|^2 dx \left|_{t=0}^T \right. + c(\eta,T)\varepsilon^{\frac{1}{4}}.
\end{equation}
From (\ref{pr1}) to (\ref{pr6}) we conclude that
$$
\frac{1}{\varepsilon^{2}}\int_0^T \int_{\mathbb{R}^3}\left(\varrho_\varepsilon-r\right)\partial_{t}H^{\prime}\left(r\right) - p\left(\varrho_\varepsilon \right)\textrm{div}\mathbf{U}_{\varepsilon,\eta} - \varrho_\varepsilon\mathbf{u}_\varepsilon\cdot \nabla H^{\prime}\left(r\right) dxdt
$$
\begin{equation}\label{prfinal}
\leq \frac{1}{2} \int_{\mathbb{R}^3} a^2 \left|\psi_{\varepsilon,\eta}\right|^2 dx \left|_{t=0}^T \right. + c(\eta,T)\varepsilon^{\alpha}, \, \alpha=\min\{\frac{1}{4},\frac{1}{2\gamma}\}.
\end{equation}

\bigskip
\bigskip
\textbf{Acknowledgements}

M. C. has been fully supported by the Croatian Science Foundation under the project MultiFM IP-2019-04-1140. The author would like to thank professor Igor Pa\v zanin and professor Boris Muha from the Department of Mathematics, Faculty of Science, University of Zagreb for the fruitful discussions about the topic.

\end{document}